\documentclass[12pt]{article}
\usepackage[english,russian]{babel}
\usepackage[cp1251]{inputenc}
\usepackage{amsfonts,amssymb,latexsym,amsmath,enumerate}
\usepackage{amsfonts}
\usepackage{indentfirst}
\usepackage{graphicx,graphics,verbatim}

\usepackage{wrapfig}
\usepackage{geometry} 
\usepackage{cite}
\usepackage{cmap}

\newtheorem{teo}{Теорема}
\newtheorem{pr}{Пример}
\newtheorem{zd}{Задача}

\begin{document}
	
		MSC 26D20; 97I40; 97I50

\vspace{2cm}

\begin{Large}	
\begin{center}		
	{\bf INEQUALITIES IN CALCULUS:\\}
\vspace{1cm}
 {\bf METHODS OF PROVING RESULTS AND PROBLEM SOLVING} 
\end{center}	
\end{Large}

\vspace{2cm}

\begin{large}
\begin{center}
\textbf{Sergei  Sitnik}\\
	Belgorod State University, Belgorod, Russia\\
	sitnik@bsu.edu.ru\\	
\vspace{1cm}
	
\textbf{Elina  Shishkina}\\
	Voronezh State University, Voronezh, Russia\\
	shishkina@amm.vsu.ru\\
\vspace{1cm}

\textbf{Lidiya Kovaleva}\\
	Belgorod State University, Belgorod, Russia\\
	kovaleva\_l@bsu.edu.ru\\	
\vspace{1cm} 

\textbf{Olga Chernova}\\
	Belgorod State University, Belgorod, Russia\\
	chernova\_olga@bsu.edu.ru	
\end{center}
\end{large}

\newpage


\vspace{1cm}

{\bf Abstract}. This preprint is a text for students and teachers on inequalities. Some standard topics are covered on application of calculus to inequality proving. Many examples are considered, stated, solved or partially solved. Some problems are standard, but some are rare, new and original. 

The next topics are considered with many examples: monotonicity of functions, Lagrange theorem and inequalities proving, estimating of finite sums, inequalities of Schl\"omilch --- LeMonnier type, proof of inequalities by  method of mathematical induction, inequalities for the number $e$, exponentials, logarithmic  and similar functions, some means and their inequalities, Cauchy--Bunyakovskii, Minkovskii, Young, H\"older (Rogers--H\"older--Riesz !) inequalities and some of their improvements and generalisations.

Some new results include inequalities on exponentials, logarithmic  and similar functions, generalisations of Cauchy--Bunyakovskii and Young inequalities, some mean inequalities including mean inequalities on the complex domain and more.

Abstracts and table of contents are in English, main text is in Russian.
\newpage

\begin{center}
 \bf{ TABLE  OF  CONTENTS}
\end{center}

\begin{enumerate}

\item Introduction ---           p. 6.

\item Monotonicity of functions, Lagrange theorem and inequalities proving ---           p. 7.

\item Estimating of finite sums  --- p. 16.

\item Inequalities of Schl\"omilch --- LeMonnier type  p. 18.

\item Proof of inequalities by  method of mathematical induction --- p. 26.

\item Inequalities for the number $e$, exponentials  and similar functions --- p. 29.

\item Means and their inequalities --- p. 41. 

\begin{enumerate}

\item Different kinds of mean values --- p. 41. 

\item Inequalities for means in the complex domain --- p. 45.

\item Inequalities for means --- p. 48.

\end{enumerate}

\item Cauchy--Bunyakovskii, Minkovskii, Young, H\"older (Rogers--H\"older--Riesz !) inequalities and some of their improvements and generalisations --- p. 49.
    
\item References --- p. 54.

\end{enumerate}
\newpage

\begin{Large}	
\begin{center}		
	{\bf НЕРАВЕНСТВА В АНАЛИЗЕ:\\}
\vspace{1cm}
 {\bf методы доказательства результатов и решения задач}
\end{center}	
\end{Large}

\vspace{2cm}

\begin{large}
\begin{center}
\textbf{Сергей Ситник}\\
	Белгородский государственный национальный исследовательский университет (НИУ "БелГУ"), Белгород, Россия\\
	sitnik@bsu.edu.ru\\	
\vspace{1cm}
	
\textbf{Элина Шишкина}\\
	Воронежский государственный университет, Воронеж, Россия\\
	shishkina@amm.vsu.ru\\
\vspace{1cm}

\textbf{Лидия Ковалёва}\\
	Белгородский государственный национальный исследовательский университет (НИУ "БелГУ"), Белгород, Россия\\
	kovaleva\_l@bsu.edu.ru\\	
\vspace{1cm}

\textbf{Ольга Чернова}\\
	Белгородский государственный национальный исследовательский университет (НИУ "БелГУ"), Белгород, Россия\\
	chernova\_olga@bsu.edu.ru	
\end{center}
\end{large}

\newpage

\tableofcontents
\newpage

\section{Введение }

Издание содержит основные факты, на которые опирается доказательство неравенств: теоремы о среднем, свойства монотонности функций, оценки сумм и рядов, дискретные и интегральные  классические неравенства. В него  также включены методические указания, задания для самостоятельной работы студентов, образцы решения модельных задач по данному предмету.
Текст включает в себя большое количество задач, нацеленных на закрепление и углубление знаний по математическому анализу, часть из этих задач приводится в учебной литературе впервые.

Для понимания условий задач и их решения достаточно подготовки в объёме основных сведений из курса дифференциального и интегрального исчисления Г. М. Фихтенгольца \cite{Fikh}-\cite{Fikh2}. С материалом по решению неравенств можно познакомится по книге \cite{ineq}. При подготовке этого учебно--методического пособия использовался материал из препринта \cite{Sit1} и обзора \cite{Sit2}, а также  некоторые материалы с сайта мехмата МГУ dxdy. Для знакомства с некоторыми более продвинутыми неравенствами рекомендуем обратиться к литературе в конце препринта, в том числе к  работам \cite{Sit1}, \cite{Sit2}--\cite{Sit15}.

В тексте некоторые задачи помечены звёздочками, это задачи повышенной трудности. Уровень сложности указан количеством звёздочек: * --- трудная задача, ** --- очень трудная задача, *** --- исследовательская задача, решение которой может быть неизвестно авторам. Понятно, что эта классификация носит достаточно условный характер.

Данный текст был издан в 2021 году авторами как научно--методическое пособие для студентов Белгородского государственного университета \cite{pos}.

\newpage

\section{Монотонность функции, теорема Лагранжа и доказательства неравенств}

Исследование неравенств с помощью производных основывается на следующих теоремах.

\begin{teo}\label{teo1} Пусть функции $f(x)$ и $g(x)$ определены  и дифференцируемы на $(a, b)$
	и выполненяются два условия:
\begin{equation}\label{Eq01}
	\lim\limits_{x\rightarrow a+0} f(x) \geq \lim\limits_{x\rightarrow a+0} g(x),
\end{equation}
\begin{equation}\label{Eq02}
	 f'(x) \geq g'(x),\qquad \forall x \in (a, b).
\end{equation}
Тогда справедливо неравенство
\begin{equation}\label{Eq03}
f(x) \geq g(x),\qquad \forall x \in (a, b).
\end{equation}

Если неравенство \eqref{Eq02} выполняется со знаком "$>$"\,, то и неравенство \eqref{Eq03} выполняется со знаком "$>$"\,.
\end{teo}

{\it Доказательство:} Действительно, поскольку	функции $f(x)$ и $g(x)$ дифференцируемы на $(a, b)$ и для их производных выполняется неравенство \eqref{Eq03}, то
$$
f'(x)-g'(x)\geq0,\qquad \forall x \in (a, b)
$$
а это означает, что функция $f(x)-g(x)$ возрастает на $(a,b)$ и, следовательно,
$$
f(x)-g(x) \geq \lim\limits_{x\rightarrow a+0}(f(x)-g(x)) \qquad \forall x \in (a, b).
$$
В силу \eqref{Eq01} имеем $\lim\limits_{x\rightarrow a+0}(f(x)-g(x)) \geq 0$, тогда
$$
f(x)-g(x) \geq \lim\limits_{x\rightarrow a+0}(f(x)-g(x))\geq 0 \qquad \forall x \in (a, b),
$$
а это и означает, что неравенство \eqref{Eq03} выполняется.

{\it Замечание 1.}
Поскольку из условия $f'(x)-g'(x)>0$ следует строгое возрастание функции $f(x)-g(x)$  на $(a,b)$
$$
f(x)-g(x) > \lim\limits_{x\rightarrow a+0}(f(x)-g(x)) \qquad \forall x \in (a, b),
$$
то в этом случае неравенство \eqref{Eq03} выполняется со строгим знаком "$>$"\,.

\vskip 0.5 cm

{\it Замечание 2.} Если в теореме \ref{teo1} функции $f(x)$ и $g(x)$ определены в точке $a$, то вместо условия \eqref{Eq01} необходимо выполнение  условия
$$
	 f(a) \geq  g(a).
$$

\vskip 0.5 cm

{\it Замечание 3.} В теореме \ref{teo1} допустимо рассматривать и $b=+\infty$.

\vskip 0.5 cm

Из условия монотонности следует и более простое, но удобное утверждение.

\begin{teo}\label{teo2} Пусть функция $h(x)$  дифференцируема на $(a,b)$. Тогда, если $h'(x) \geq 0, \forall x \in (a, b)$, то
\begin{equation}\label{Eq04}
\lim_{x\to a+0}h(x) \leq h(x) \leq \lim_{x\to b-0}h(x).
\end{equation}
Если справедливо строгое неравенство $h'(x)>0$ для всех $x \in (a, b)$, то и \eqref{Eq04} выполняется со строгими знаками неравенств при $x\in(a,b)$.
\end{teo}

Рассмотрим примеры доказательства неравенств при помощи теорем \ref{teo1} и \ref{teo2}.
\begin{pr}\label{pr0}
	Доказать неравенство
	$$
	\sqrt{1+x}\leq 1+\frac{x}{2}
	$$
	при $x\geq 0$.
\end{pr}
{\it Решение.} Пусть  $f(x)=1+\frac{x}{2}$, $g(x)=\sqrt{1+x}$. Тогда неравенство \eqref{Eq01} выполняется:
$$
	 f(0)= 1\geq 1 = g(0).
$$
Поскольку $f'(x)=\frac{1}{2}$, $g'(x)=\frac{1}{2\sqrt{1+x}}$, то неравенство \eqref{Eq02} также выполняется:
$$
	f'(x)=\frac{1}{2} \geq \frac{1}{2\sqrt{1+x}}=g'(x),\qquad \forall x\geq 0.
$$
Следовательно, по теореме \ref{teo1} справедливо неравенство
$$
f(x)=1+\frac{x}{2} \geq \sqrt{1+x}=g(x),\qquad \forall x\geq 0.
$$

\begin{pr}\label{pr01}
Определить, что больше $e^\pi$ или $\pi^e$.
\end{pr}
{\it Решение.}
Рассмотрим функцию $h(x)=-\frac{\ln x}{x}$ на интервале $(e,\pi)$. Найдем её производную:
$$
h'(x)=\left(-\frac{\ln x}{x} \right)'=-\frac{\frac{1}{x}\cdot x-\ln x}{x^2}=\frac{\ln x-1}{x^2}>0,\qquad \forall x>e.
$$
Тогда по теореме \ref{teo2} имеем
$$
\lim_{x\to e+0}h(x)=-\frac{\ln e}{e} < -\frac{\ln x}{x}=h(x) < -\frac{\ln \pi}{\pi}=\lim_{x\to \pi-0}h(x).
$$
Следовательно,
$$
\frac{\ln e}{e}>\frac{\ln \pi}{\pi}\,\Rightarrow \, \pi \ln e>e\ln \pi \,\Rightarrow \,  \ln e^\pi>\ln \pi^e \,\Rightarrow \, e^\pi>\pi^e.
$$
Таким образом, $e^\pi$ больше, чем $\pi^e$.

\begin{zd}\label{zd01}
Доказать, что при любом действительном $x$ справедливо неравенство
$$
|\sin{x}|\leq|x|.
$$
\end{zd}

\begin{zd}\label{zd02}
	Определить, что больше $100^{101}$ или $101^{100}$.
\end{zd}

\begin{zd}\label{zd03}
Доказать, что для любых положительных чисел
$a$ и $b$
справедливо неравенство
$$
\frac{a}{b}+\frac{b}{a}\geq2.
$$
\end{zd}
{\it Указание:} введите обозначение $\frac{a}{b}=x\geq 1$ и примените теорему \ref{teo1} к функциям $f(x)=x$ и $g(x)=2-\frac{1}{x}$.

\begin{zd}\label{zd04}
Доказать, что при $x>0$ справедливо неравенство
$$
1+\ln x^2\leq x^2.
$$
\end{zd}

\begin{zd}\label{zd05}
Показать, что при $x>1$
выполняется неравенство
$$
\sqrt{x}+\frac{1}{\sqrt{x}}>2.
$$
\end{zd}

\begin{zd}\label{zd06}
Доказать, что при
$x\neq 0$
выполняется неравенство
$$
e^x>1+x.
$$
\end{zd}

\begin{zd}\label{zd07}
	Доказать, что на интервале $\left(0,\frac{\pi}{2} \right) $
		выполняется неравенство
	$$
	\sin x+\tg x> x.
	$$
\end{zd}

\vskip 0.5 cm
Приведем примеры, в которых эта теорема используется для доказательства логарифмических неравенств.
\vskip 0.5 cm

\begin{pr}\label{pr1}
Доказать, что справедливо неравенство
\begin{equation}\label{Eq05}
\ln(1+x) \leq x, \quad x \geq 0,
\end{equation}
причем знак равенства возможен лишь при $x = 0$.
\end{pr}
{\it Решение.} Введем обозначения
$$
f(x) = x,\qquad g(x) = \ln(1+x),\qquad a = 0.
$$

Очевидно, что при $x = a = 0$
$$
f(0) = g(0) = 0.
$$

Следовательно, выполнено \eqref{Eq01}. Поскольку
$$
f'(x) = 1, \qquad g'(x) = \frac{1}{1+x},
$$
то неравенство для производных \eqref{Eq02} принимает вид
$$
1 > \frac{1}{1+x}, \quad x>0.
$$

Это неравенство выполняется со строгим знаком "$>$"\, при всех $x>0$, а равенство возможно только при $x=0$. Следовательно, по Tеореме \ref{teo1} неравенство \eqref{Eq05} доказано.

\begin{pr}\label{pr2} Доказать неравенство
\begin{equation}\label{Eq06}
\frac{x}{x+1} \leq \ln(1+x),\qquad x\geq0,
\end{equation}
причем знак равенства возможен лишь при $x=0$.
\end{pr}
{\it Решение.}  Аналогично предыдущему примеру обозначаем
$$
f(x) = \ln(1+x),\qquad g(x) = \frac{x}{x+1}.
$$

Так как $f(0)=0\geq0=g(0)$, то осталось проверить неравенство \eqref{Eq02}. Имеем
$$
f'(x) = \frac{1}{1+x}, \qquad g'(x) = \frac{1}{(1+x)^2},
$$
и, очевидно,
$$
\frac{1}{1+x} > \frac{1}{(1+x)^2},\qquad x>0.
$$
Равенство $\frac{1}{1+x}=\frac{1}{(1+x)^2}$ возможно только при $x=0$.
Следовательно, по  Tеореме \ref{teo1} неравенство \eqref{Eq06} доказано.

\begin{pr}\label{pr3} Доказать, что при $x>0$ справедливо неравенство
\begin{equation}\label{Eq07}
\ln^2\left(1+\frac{1}{x}\right)<\frac{1}{x(x+1)}.
\end{equation}
\end{pr}
{\it Решение.} Здесь затруднительно применить теорему \ref{teo1}. Поэтому введем функцию
$$
h(x) = x \cdot \ln^2\left( 1 + \frac{1}{x}\right) - \frac{1}{x+1}.
$$
Неравенство \eqref{Eq07} эквивалентно условию $h(x) < 0$ при $x>0$. Вычисляем производную
$$
h'(x) = \ln^2\left( 1 + \frac{1}{x}\right)  + x \cdot 2 \cdot \ln\left( 1 + \frac{1}{x}\right)  \cdot \frac{1}{1 + \frac{1}{x}} \cdot \left( - \frac{1}{x^2}\right)  + \frac{1}{(x+1)^2} =
$$
$$
= \ln^2\left( 1 + \frac{1}{x}\right)  - 2 \cdot \ln\left( 1 + \frac{1}{x}\right)  \cdot \frac{1}{1+x} + \frac{1}{(x+1)^2} =
$$
$$
=\bigg[\ln\left( 1 + \frac{1}{x}\right)  - \frac{1}{x + 1}\bigg] \geq 0.
$$

Далее мы покажем, что при $x>0$ справедливо
$$
\ln\left( 1 + \frac{1}{x}\right)  > \frac{1}{x+1},
$$
поэтому $h'(x)>0$ строго при $x>0$. Следовательно, по теореме \ref{teo2}
$$
h(x) < \lim_{x\to +\infty}\bigg[x \cdot \ln^2\left(1 + \frac{1}{x}\right) - \frac{1}{x+1}\bigg]=
$$
$$
 =\lim_{x\to +\infty}\bigg[\ln\left( 1+\frac{1}{x}\right) ^{x} \cdot \ln\left(1+\frac{1}{x}\right) - \frac{1}{x+1}\bigg]
= 1 \cdot 0 - 0 = 0.
$$
Итак, $h(x) < 0$ откуда и  следует \eqref{Eq07}.

\begin{zd}\label{zd1}
 Доказать, что справедливо двустороннее неравенство
\begin{equation}\label{Eq08}
\displaystyle -1 \leq x \ln^2\left( 1+\frac{1}{x}\right) - \frac{1}{x+1} < 0, \quad x\geq 0,
\end{equation}
причем равенство в левой части достигается лишь при $x \rightarrow 0$.
\end{zd}
\vskip 0.5cm

Константы $-1$ и $0$ в неравенстве \eqref{Eq08} являются точными. Это означает, что неравенство станет неверным для всех $x>0$, если $-1$ заменить любым большим числом, а $0$ --- любым меньшим.

Сделаем полезное замечание. Если попытаться доказать \eqref{Eq07} в виде
$$
\frac{1}{x+1} > x \cdot \ln^2\left( 1 + \frac{1}{x}\right) , \quad x > 0,
$$
по теореме \eqref{teo1}, то условие \eqref{Eq01} выполнено, а условие \eqref{Eq02} принимает вид
$$
- \frac{1}{(x+1)^2}>\ln^2\left( 1+\frac{1}{x}\right)  + x \cdot 2 \cdot \ln\left( 1+\frac{1}{x}\right)  \cdot \frac{1}{1 + \frac{1}{x}} \cdot \left( -\frac{1}{x^2}\right)  =
$$
$$
= \ln^2\left( 1 + \frac{1}{x}\right)  - 2 \cdot \ln\left( 1+\frac{1}{x}\right)  \cdot \frac{1}{x+1} \qquad \Rightarrow
$$
$$
 \bigg[\ln\left( 1+\frac{1}{x}\right)  - \frac{1}{x+1}\bigg]^2 < 0.
$$
Но последнее неравенство не выполняется ни при каких $x$.

Разумеется, это не означает, что мы доказываем неверное неравенство \eqref{Eq07}. Просто доказательство невозможно тем методом, который мы выбрали.

Интересно все--таки поискать такое доказательство неравенства \eqref{Eq07}, которое использует теорему \ref{teo1}.
Для этого перепишем \eqref{Eq07} в виде
\begin{equation}\label{Eq081}
\ln(1+y) \leq \frac{y}{\sqrt{1+y}}, \quad y = \frac{1}{x} > 0.
\end{equation}

\begin{pr}\label{pr4} (второе доказательство неравенства \eqref{Eq07}).
Доказать, что неравенство \eqref{Eq081} справедливо, причем равенство достигается лишь при $y=0$.
\end{pr}
{\it Решение}. Обозначим
$$
g(y) = \ln(1+y), \qquad f(y) = \frac{y}{\sqrt{1+y}},\qquad a = 0.
$$

Очевидно, что условие \eqref{Eq01} выполнено, и осталось сравнить производные
$$
\frac{1}{1+y} \leq \frac{\sqrt{1+y} - \frac{y}{2\sqrt{1+y}}}{1+y} = \frac{y+2}{2(1+y)\sqrt{1+y}} \Leftrightarrow
$$
$$
\Leftrightarrow 2\sqrt{1+y} \leq y+2 \Leftrightarrow \sqrt{1+y} \leq 1 + \frac{y}{2} \Leftrightarrow y^2 \geq 0.
$$

Равенство в последнем случае достигается лишь при $y = 0$, что и требовалось доказать.

\vskip 0.5 cm

Отметим, что соотношение
$$
\sqrt{1+y} \leq 1 + \frac{y}{2}, \quad y > 0,
$$

есть частный случай неравенства Бернулли.

Исследуем, является ли константа $1$ в правой части \eqref{Eq081} точной. Иными словами, не будет ли справедливо более точное неравенство с уменьшенной по сравнению с \eqref{Eq081} правой частью
\begin{equation}\label{Eq09}
\ln(1+y) \leq \frac{y}{y+\delta}, \quad y \geq 0,
\end{equation}
где $\delta$ некоторая постоянная, $\delta>1$?

\begin{pr}\label{pr5} Доказать, что постоянная $1$ в правой части \eqref{Eq081} является точной, то есть неравенство \eqref{Eq09} не может выполнять для всех $y \geq 0$ ни при каком $\delta>1$.
\end{pr}
{\it Решение.} При доказательстве используется очень полезный технический прием, которым надо овладеть. Этот метод введения вспомогательной функции. Таким приемом многие неравенства можно тестировать на точность и находить пути их улучшения.

Определим функцию $\varepsilon(y)$ равенством
$$
\ln(1+y) = - \frac{y}{\sqrt{y+\varepsilon(y)}},\qquad \varepsilon(y) = \frac{y^2}{\ln^2(1+y)} - y.
$$

Найдем предел
$$
\lim\limits_{y\to 0} \varepsilon (y)= \lim\limits_{y\to 0}\bigg(\dfrac{y^2}{\ln^2 (1+y)}-y\bigg)=1-0=1
$$
Поэтому, неравенство \eqref{Eq09} невозможно.
Действительно, предположим противное, что выполнено \eqref{Eq09} при некотором $\delta >1$. Тогда получаем

$$
\ln (1+y)=\dfrac{y}{\sqrt{y+\varepsilon(y)}}\leq \dfrac{y}{\sqrt{y+\delta}} \Rightarrow \varepsilon(y) \geq \delta>1,
$$
но это невозможно хотя бы при $y \to 0$. Полученное противоречие завершает решение.

\vskip 0.5 cm

\begin{zd}\label{zd2} Доказать, что
$$
\lim\limits_{y \to +\infty}\varepsilon (y)=\lim\limits_{y \to +\infty}\left(\frac{y^2}{\ln^2(1+y)} - y\right) = +\infty.
$$
\end{zd}

\begin{zd}\label{zd3} Докажите, что функция $\varepsilon(y)$ монотонно возрастает при $y>0$. Постройте её график с использованием одного из пакетов компьютерной математики.

\begin{zd}\label{zd5} Доказать, что при $x>0$
\begin{equation}\label{Eq10}
\dfrac{2x}{x+2}<\ln(1+x).
\end{equation}
\end{zd}

\begin{zd}\label{zd6} Доказать, что постоянная 2 в знаменателе \eqref{Eq10} точная.
\end{zd}
{\it Указание:} введите вспомогательную функцию $\varepsilon=\varepsilon(x)$ по формуле
\begin{equation}\label{Eq11}
\dfrac{2x}{x+\varepsilon(x)}=\ln(1+x).
\end{equation}
\end{zd}

\begin{zd}\label{zd7} Доказать, что функция $\varepsilon(x)$, определяемая из \eqref{Eq11}, монотонно убывает при $x>0$, причем
$$
\lim\limits_{x \to +\infty} \varepsilon (x)=-\infty.
$$
\end{zd}
\vskip 0.5 cm

Информацию о вспомогательной $\varepsilon$-функции, как правило, всегда можно использовать для уточнения неравенства.

\begin{pr}\label{pr6}
Доказать, что при $x>a\geq0$ справедливо неравенство
$$
\ln(1+x)>\dfrac{2x}{x+\varepsilon(a)}, \qquad x>a.
$$
\end{pr}
{\it Решение.} Из определения следует, что $x+\varepsilon(x)>0$ при $x>0$. Следовательно, при $x>a$ следующие неравенства эквивалентны:
$$
\varepsilon(x)<\varepsilon(a), \qquad x>a, \qquad \ln(1+x)=\dfrac{2x}{x+\varepsilon(x)}>\dfrac{2x}{x+\varepsilon(a)},
$$
 а монотонность $\varepsilon(x)$ следует из Задачи \ref{zd7}.

В частности, например,
$$
\ln (1+x)>\dfrac{2x}{x+b}, \qquad x>1, \qquad b=\dfrac{2}{\ln2}-1.
$$

\begin{zd}\label{zd8} Доказать, что при $x>0$
\begin{equation}\label{Eq12}
x-\dfrac{x^2}{2}<\ln(1+x)<x-\dfrac{x^2}{2}+\dfrac{x^3}{3}.
\end{equation}
\end{zd}

\begin{zd}\label{zd9} Из какой теоремы Анализа следуют неравенства \eqref{Eq12} и подобные им при $0<x<1$?
\end{zd}

Интересно, что \eqref{Eq12} справедливы при всех $x>0$. Это больше того, что утверждается в Задаче \ref{zd9}.

\begin{zd}\label{zd10} Доказать, что при $x>0$ справедливо неравенство $$\ln x \leq \dfrac{x^2-1}{2}.$$
\end{zd}
{\it Указание:} рассмотрите отдельно промежутки $x \in (0,1]$ и $x \in (1,\infty)$.

\begin{zd}\label{zd11} Доказать, что при $x\geq 0$ справедливо неравенство
\begin{equation}\label{Eq13}
\ln (1+x) \leq \dfrac{x(x+2)}{2(x+1)},
\end{equation}
причем равенство достигается лишь при $x=0$.
\end{zd}

\begin{zd}\label{zd12} Доказать, что постоянные $1$ и $2$ в выражениях $x+2$, $x+1$ в правой части неравенства \eqref{Eq13} --- точные.
\end{zd}

\begin{zd}\label{zd13} Докажите, основываясь на неравенствах \eqref{Eq081} и \eqref{Eq13}, что справедлива цепочка неравенств
\begin{equation}\label{Eq14}
\ln (1+x)\leq\dfrac{x}{\sqrt{1+x}}\leq\dfrac{x(x+2)}{2(x+1)}\leq x, \quad x\geq 0.
\end{equation}
\end{zd}
Таким образом, неравенство \eqref{Eq081} самое сильное из доказанных, а неравенство \eqref{Eq13} сильнее \eqref{Eq05}.

\begin{zd}\label{zd14}
Докажите, что полученные неравенства \eqref{Eq06}  и \eqref{Eq10} можно расположить в последовательность
\begin{equation}\label{Eq15}
\dfrac{x}{x+1}\leq \dfrac{2x}{x+2}\leq \ln (1+x) , \qquad x\geq0.
\end{equation}
\end{zd}

\begin{zd}\label{zd15} Сравнить, какие из неравенств \eqref{Eq12}, \eqref{Eq14}, \eqref{Eq15}  сильнее.
\end{zd}

Примеры 3--8 и Задачи 8--13, 17--20 относятся к теме об оценках логарифма через рациональные функции. Такие оценки систематически получаются при помощи так называемых аппроксимаций Паде, то есть рациональных аппроксимаций функции, которые имеют наилучший порядок приближения в нуле или на бесконечности.

\begin{zd}\hspace{-2mm}***\label{zd16}
 Изучить неравенства \eqref{Eq14}--\eqref{Eq15} между модулями функций в комплексной плоскости.
\end{zd}

\begin{zd}\hspace{-2mm}***\label{zd17} При каких $z \in \mathbb{C}$ справедливо неравенство
\begin{equation}\label{Eq16}
|\ln (1+z)|\leq|z|?
\end{equation}
Здесь $\ln (1+z)$ --- главное значение комплексного логарифма,  то есть $\ln (1+z)=\ln r\,\cdot \exp(i\varphi)$, где $r=|z+1|, \, \varphi={\rm arg}(z+1), \varphi \in (-\pi,\pi]$.
\end{zd}

Авторам неизвестно решение Задачи \ref{zd17}.

{\it Гипотеза:} на комплексной плоскости существует замкнутая кривая, на которой в \eqref{Eq16} достигается равенство, при этом внутри кривой неравенство \eqref{Eq16} выполняется со  знаком "$>$"\,, а вне кривой --- со знаком "$<$"\,.

\begin{zd}\label{zd18}
Постройте на компьютере все такие $z \in \mathbb{C}$, для которых в \eqref{Eq16} достигается равенство.
\end{zd}

\begin{zd}\label{zd19}
 Докажите, что неравенство \eqref{Eq16} неверно $\forall z \in \mathbb{C}$.
\end{zd}
{\it Указание:} рассмотрите $z$ на прямой $z=-1+iy$, $y\in \mathbb{R}$.

\vskip 0.5 cm

 Перепишем неравенства\eqref{Eq12}, \eqref{Eq14}--\eqref{Eq15}, произведя  замену переменной $x=\dfrac{1}{y}$. Получим
\begin{equation}\label{Eq17}
\ln \left( 1+\dfrac{1}{y}\right)  < \dfrac{1}{\sqrt{y(y+1)}} < \dfrac{2y+1}{2y(y+1)}<\dfrac{1}{y}, \qquad y>0,
\end{equation}
\begin{equation}\label{Eq18}
\ln \left( 1+\dfrac{1}{y}\right) < \dfrac{1}{y+\frac{1}{2}} < \dfrac{1}{y+1}, \qquad  y>0,
\end{equation}
\begin{equation}\label{Eq19}
\dfrac{1}{y} -\dfrac{1}{2y^2} < \ln\left( 1+\dfrac{1}{y}\right) <\dfrac{1}{y}-\dfrac{1}{2y^2}+\dfrac{1}{3y^3}.
\end{equation}
Приведенные неравенства \eqref{Eq17}--\eqref{Eq19} часто используются в Математическом анализе.

В частности, из  неравенств \eqref{Eq17}--\eqref{Eq18} следует, что
\begin{equation}\label{Eq20}
\dfrac{1}{y+1}<\ln\left( 1+\dfrac{1}{y}\right) <\dfrac{1}{y}, \qquad y>0.
\end{equation}

\vskip 0.5 cm

Хорошо известная теорема Лагранжа о среднем (или о промежуточном значении) --- одна из самых важных теорем математического анализа. В этом разделе покажем, как применять ее к доказательству неравенств.

\begin{teo} (Лагранж).
	Если функция $f$ непрерывна на отрезке $[a,b]$ и в каждой точке интервала $(a,b)$ имеет производную,
	 то в этом интервале существует по крайней мере одна такая точка $\xi$, что
	\begin{equation}\label{Lag1}
	f'(\xi)=\frac{f(b)-f(a)}{b-a}.
\end{equation}
\end{teo}

Положим $a=x$, $b-a=\Delta x$ и $b=x+\Delta x$, тогда теорему  Лагранжа можно переформулировать так:
\begin{equation}\label{Lag2}
f(x+\Delta x)-f(x)=f'(x+\theta \Delta x)\Delta x,\qquad 0<\theta<1.
\end{equation}
Эта формула называется \textbf{формулой конечных приращений Лагранжа}.

\begin{pr}\label{Lag01}
Доказать неравенство
$$
|{\rm arctg}\,a-{\rm arctg}\,b |\leq|a-b|.
$$
\end{pr}
{\it Решение.} Рассмотрим функцию
$f(x)={\rm arctg}\,x$ на отрезке $[a,b]$. Для этой функции справедлива теорема Лагранжа и $f'(x)=\frac{1}{1+x^2}$.
Следовательно, применяя теорему Лагранжа, получим
$$
\frac{{\rm arctg}\,b-{\rm arctg}\,a}{b-a}=\frac{1}{1+\xi^2},\qquad \xi\in(a,b)
$$
или
$$
{\rm arctg}\,a-{\rm arctg}\,b=(a-b)\frac{1}{1+\xi^2}.
$$
Перейдем в данном равенстве к модулям:
$$
|{\rm arctg}\,a-{\rm arctg}\,b|=|a-b|\frac{1}{1+\xi^2}.
$$
Поскольку $\frac{1}{1+\xi^2}\leq 1$, то
$$
|{\rm arctg}\,a-{\rm arctg}\,b|\leq|a-b|.
$$

\begin{zd}\label{zdLag1}
	Доказать неравенство
$$
|\sin a-\sin b|\leq|a-b|.
$$
\end{zd}

\begin{zd}\label{zdLag2}
	При $0<a\leq b$ доказать  неравенства
	$$
	\frac{b-a}{b}\leq\ln\frac{b}{a}\leq\frac{b-a}{a}.
	$$
\end{zd}

\section{Оценки конечных сумм}

В этом разделе будем рассматривать  оценки сумм вида
$$
\sum\limits_{k=1}^n a_k,
$$
где $a_k>0$, $\lim\limits_{k\rightarrow+\infty}a_k=0$.

\begin{pr}\label{pr7}
Докажите, что если обозначить
\begin{equation}\label{Eq21}
\displaystyle S_n=1+\dfrac{1}{2}+\dfrac{1}{3}+...+\dfrac{1}{n}=\sum_{k=1}^n \dfrac{1}{k},
\end{equation}
то для любого $n \in \mathbb{N}$ будут выполнены неравенства
\begin{equation}\label{Eq22}
\ln(1+n)<S_n<\ln(n+1)+\dfrac{n}{n+1}.
\end{equation}
\end{pr}
{\it Решение.} Запишем неравенства \eqref{Eq20} для всех $x=k$, $1\leq k \leq n$:
$$
\dfrac{1}{2}<\ln\left( 1+\dfrac{1}{1}\right) <\dfrac{1}{1},
$$
$$
\dfrac{1}{3}<\ln\left( 1+\dfrac{1}{2}\right) <\dfrac{1}{2},
$$
$$
\dfrac{1}{k+1}<\ln\left( 1+\dfrac{1}{k}\right) <\dfrac{1}{k},
$$
$$
\dfrac{1}{n+1}<\ln\left( 1+\dfrac{1}{n}\right) <\dfrac{1}{n}.
$$
Складывая все эти неравенства, получаем
$$
S_n -1+\dfrac{1}{n+1}<\sum_{k=1}^{n}\ln\left( 1+\dfrac{1}{k}\right)<S_n.
$$
Но сумму логарифмов можно вычислить
$$
\displaystyle \sum_{k=1}^{n}\ln\left( 1+\dfrac{1}{k}\right)=\sum_{k=1}^{n}\big[\ln(k+1)-\ln(k)\big]=(\ln2-\ln1)+(\ln3-\ln2)+\ldots
$$
$$
\ldots+(\ln(n+1)-\ln(n))=\ln(n+1).
$$
Поэтому
$$
S_n<\ln(n+1)+1-\dfrac{1}{n+1}=\ln(n+1)+\dfrac{n}{n+1}, \quad S_n>\ln(n+1),
$$
что и требовалось доказать.

Приём, применённый при доказательстве этого неравенства, называется "телескопическим методом"{}, он заключается в сокращении основной части слагаемых при суммировании, аналогично тому, как складывается труба телескопа, уменьшая свою длину.

\vskip 0.5 cm

Одним из следствий неравенства \eqref{Eq22} является то, что сумма $S_n$ стремится к бесконечности при $n \to +\infty$. Действительно,
$$
+\infty>S_n>\ln(1+n)=\lim\limits_{n \to +\infty} S_n=+\infty
$$
по лемме "о двух милиционерах"\  (они же два полицейских или встречные последовательности)\, для расширенной числовой прямой. Сумма $S_n$, определяемая по формуле \eqref{Eq21}, является частичной суммой гармонического ряда
\begin{equation}\label{GR}
\sum\limits_{n=1}^\infty \frac{1}{n}.
\end{equation}
Таким образом, мы показали, что гармонический ряд, \eqref{GR} расходится, то есть
\begin{equation}\label{GRRas}
\lim\limits_{n \to +\infty}S_n{=}{+}\infty.
\end{equation}

\begin{zd}\label{zd20}
Докажите, что существует предел
\begin{equation}\label{EqC}
\lim \limits_{n\to+\infty}(S_n-\ln(n))=C, \qquad 0<C<1.
\end{equation}
\end{zd}
{\it Указание:} последовательность $S_n -\ln(n)$ монотонно убывает.

Из формулы Задачи \ref{zd20} следует, что при больших $n$ для величины $S_n$ справедливо приближённое асимптотическое равенство $S_n{~\approx~}\ln(n){+}C$, и, следовательно, мы вновь получаем $S_n \to \infty$. Константа $C$ называется постоянной Эйлера. Можно вычислить, что $C\approx 0,577$. До сих пор не решена задача о том, является ли $C$ рациональным или иррациональным числом. Эту задачу не смогли решить Эйлер и Гаусс, не поддается она и усилиям современных математиков.

\section{Неравенства типа  Шлёмильха--Лемоннье}

\begin{zd}\label{zd21}
Докажите двустороннее неравенство
\begin{equation}\label{Eq23}
\ln(n+1)<S_n<1+\ln(n+1).
\end{equation}
\end{zd}

\vskip 0.5 cm

Неравенство \eqref{Eq23} исторически было первым на заданную тему. Его предложил в качестве задачи в 1885 г. Оскар Шлёмильх (автор одной из формул для остатка ряда Тэйлора). Неравенство \eqref{Eq23} также доказал в 1889 Г. Лемоннье вместе со следующим уточнением.

\begin{zd}\label{zd22}
 Справедливо двустороннее неравенство
\begin{equation}\label{Eq24}
\ln(n+1)<S_n\leq 1+\ln n,
\end{equation}
причем равенство в правой части возможно лишь при $n=1$.
\end{zd}

Рассмотрим несколько задач на тему неравенства Шлёмильха--Лемоннье и его уточнений.

\begin{zd}\label{zd23}
Докажите двустороннее неравенство
\begin{equation}\label{Eq25}
C+\ln n<S_n\leq 1+\ln n.
\end{equation}
Постоянная $C$ определяется равенством \eqref{EqC}. Равенство в правой части \eqref{Eq25} достигается лишь при $n=1$. Постоянные $C$ и $1$ в обеих частях являются точными.
\end{zd}

\begin{zd}\label{zd24}
Докажите двустороннее неравенство
\begin{equation}\label{Eq26}
1-\ln2+\ln(n+1)\leq S_n<C+\ln(n+1).
\end{equation}
Равенство в левой части достигается лишь при $n=1$. Постоянные в обеих частях являются точными.
\end{zd}

Отметим, что $1-\ln2\approx 0,30669, \, C\approx 0,5772$.

\begin{zd}\label{zd25}
Докажите, что при $n\geq 4$ справедливо двустороннее неравенство
\begin{equation}\label{Eq27}
C+\ln n < S_n < C+\ln(n+1).
\end{equation}
\end{zd}

Неравенства \eqref{Eq27} --- это гибрид неравенств  \eqref{Eq25} и \eqref{Eq26}, вобравший в себя наиболее точную информацию из них и неулучшаемый в приведенных терминах.

\vskip 0.5 cm

Теперь рассмотрим несколько неравенств в терминах средних значений. Определения основных средних и неравенства для них будут даны далее.

Пусть
$$
M_1(a,b)=\frac{a+b}{2}, \ \ M_0(a,b)=\sqrt{ab}
$$ --- среднее арифметическое и среднее геометрическое двух неотрицательных чисел.
Тогда перепишем неравенство \eqref{Eq17} в виде
\begin{equation}\label{Eq30}
\ln\bigg(1+\frac{1}{n}\bigg) < \frac{1}{2}\bigg(\frac{1}{n} + \frac{1}{n+1}\bigg) = M_1 \bigg(\frac{1}{n},\frac{1}{n+1}\bigg).
\end{equation}

\begin{zd}\label{zd42}
Докажите неравенство
\begin{equation}\label{Eq31}
\ln(n+1) + \frac{1}{2}\frac{n}{(n+1)} < S_n,
\end{equation}
\end{zd}
\vskip 0.5 cm
Отметим, что \eqref{Eq31} лучше \eqref{Eq23}--\eqref{Eq24}, но хуже \eqref{Eq26}.

\begin{zd}\label{zd43}
Докажите неравенство
\begin{equation}\label{Eq32}
M_1\bigg(\frac{1}{x+1},\frac{1}{x+2}\bigg) < \ln\bigg(1+\frac{1}{x}\bigg), \qquad x>0.
\end{equation}
\end{zd}
\vskip 0.5 cm

Для понимания неравенств \eqref{Eq30} и  \eqref{Eq32} нужна следующая задача.

\begin{zd}\hspace{-2mm}***\label{zd44}
Исследуйте функцию  $\varepsilon(x)$, определяюмую выражением
$$
\ln\bigg(1+\frac{1}{x}\bigg) = \frac{1}{2}\bigg(\frac{1}{x+\frac{1}{2}+\varepsilon(x)}+\frac{1}{x+\frac{1}{2}-\varepsilon(x)}\bigg)=
$$
$$
=\frac{1}{2} \frac{2x+1}{\left( x+\frac{1}{2}\right)^2-\varepsilon^2(x)}.
$$
\end{zd}

\begin{zd}\hspace{-2mm}***\label{zd45}
 Можно ли усилить неравенства \eqref{Eq30}, \eqref{Eq32}, \eqref{Eq29}, заменив средние $M_0$ и $M_1$ на другие средние (см. раздел \ref{sre})?
\end{zd}

\begin{zd}\hspace{-2mm}*\label{zd46} Верно ли неравенство
$$
\ln\bigg(1+\frac{1}{x}\bigg) < M_1\bigg(\frac{1}{n},\frac{1}{n+1},\frac{1}{n+2}\bigg)?
$$
\end{zd}

\begin{zd}\hspace{-2mm}*\label{zd47}
Исследуйте последнее неравенство  с другими средними.
\end{zd}

Решения задач \ref{zd44}--\ref{zd47}  авторам неизвестны, хотя неравенство из задачи \ref{zd46} исследовать нетрудно.

\vskip 0.5 cm

Оригинальное неравенство Шлёмильха--Лемоннье было получено для гармонического ряда. Но неравенства подобного вида можно рассматривать и для других рядов, см. далее. Мы перейдем к изучению  суммы вида
$$
P_n=\sum_{k=1}^n \frac{1}{\sqrt{k(k+1)}}.
$$
При этом опять понадобятся некоторые неравенства для логарифмических функций.

\begin{zd}\label{zd26}
Докажите неравенство
\begin{equation}\label{Eq28}
\ln\bigg(1+\frac{1}{x}\bigg)>\frac{1}{\sqrt{(x+1)(x+2)}}, \qquad x>0.
\end{equation}
\end{zd}
{\it Указание:} проверьте, что идею из Примера \ref{pr3} реализовать не удается. Однако, можно рассуждать как в Примере \ref{pr4}.

\begin{zd}\label{zd27}
	Докажите, что неравенство \eqref{Eq28} самое слабое из неравенств \eqref{Eq18}.
\end{zd}

Итак, из неравенств \eqref{Eq07} и \eqref{Eq28} получаем двустореннее неравенство
\begin{equation}\label{Eq29}
\frac{1}{\sqrt{(x+1)(x+2)}} < \ln\bigg(1+\frac{1}{x}\bigg) < \frac{1}{\sqrt{x(x+1)}}, \qquad x>0.
\end{equation}

Из рассмотренных задач следует, что разные части \eqref{Eq29} не равноценны. Однако, они вместе позволяют оценить сумму $P_n$.

\begin{zd}\label{zd28}
Докажите двойное неравенство
$$
\ln(n+1) < P_n < \ln(n+1) + \frac{1}{\sqrt{2}} - \frac{1}{\sqrt{(n+1)(n+2)}},\qquad n\in\mathbb{N}.
$$
\end{zd}

\begin{zd}\label{zd29}
Докажите, что ряд
$$
\sum_{k=1}^{\infty}\frac{1}{\sqrt{k(k+1)}}
$$
расходится, то есть $\lim\limits_{n\rightarrow+\infty} P_n=\infty$.
\end{zd}

Получим теперь неравенство типа Шлёмильха--Лемоннье для $P_n$.

\begin{zd}\label{zd30}
Докажите, что последовательность
$$
a_n=P_n-\ln(n+1)
$$
монотонно возрастает.
\end{zd}

\begin{zd}\label{zd31}
 Докажите, что существует
$$
\lim_{n\rightarrow +\infty} a_n = C_1.
$$
\end{zd}

\begin{zd}\label{zd32}
Докажите, что последовательность
$$
b_n=P_n-\ln n
$$
монотонно убывает.
\end{zd}

\begin{zd}\label{zd33}
Докажите, что существует
$$
\lim_{n\rightarrow+\infty} b_n = C_2.
$$
\end{zd}

\begin{zd}\label{zd34}
Докажите равенство $C_2=C_1$.
\end{zd}

\begin{zd}\label{zd35}
Докажите двустороннее неравенство
$$
\frac{1}{\sqrt{2}} - \ln2 + \ln(n+1) \leq P_n < C_1 + \ln(n+1).
$$
Равенство в левой части достигается лишь при $n=1$.
\end{zd}

\begin{zd}\label{zd36}
Докажите двустороннее неравенство
$$
C_1 + \ln n < P_n \leq \frac{1}{\sqrt{2}} + \ln n.
$$
Равенство в правой части достигается лишь при $n=1$.
\end{zd}



\begin{zd}\label{zd39}
Доказать, что $C_1$ выражается через постоянную Эйлера по формуле
$$
C_1= C+A,
$$
где $A$ есть сумма сходящегося ряда
$$
A= \sum_{k=1}^{\infty}\frac{1}{k\sqrt{k+1}(\sqrt{k+1}+\sqrt{k})}.
$$
\end{zd}

Вернёмся к более подробному исследованию логарифмического  неравенства \eqref{Eq29}. Для анализа его точности введём вспомогательную функцию $\varepsilon(x)$.

\begin{zd}\hspace{-2mm}***\label{zd40}
Пусть функция $\varepsilon(x)$ определяется из равенства
$$
\ln^2\left(1+\frac{1}{x}\right) = \frac{1}{\sqrt{(x+1)(x+\varepsilon(x))}}.
$$
Опишите аналитически поведение $\varepsilon(x)$ на $(0,\infty)$.
Исследуйте $\varepsilon(x)$ на монотонность на  $(0,+\infty)$. Найдите пределы $\lim\limits_{x\rightarrow 0}\varepsilon(x)$ и $\lim\limits_{x\rightarrow+\infty}\varepsilon(x)$.
\end{zd}

Авторам неизвестно решение задачи \ref{zd40}.

\begin{zd}\label{zd41}
Постройте график $\varepsilon(x)$ на компьютере.
\end{zd}

\vskip 1 cm

Рассмотрим один метод получения оценок на суммы расходящегося ряда, использующий интегральный признак (см. \cite{Fikh2}, стр. 282). Пусть дан произвольный числовой знакоположительный ряд
\begin{equation}\label{Ryad}
	\sum\limits_{n=1}^{\infty} f(n),
\end{equation}
где $f(n)$ есть значение при $x=n$ некоторой функции $f(x)$, определенной для $x\geq 1$.
Функция $f(x)$ предполагается непрерывной, положительной и строго убывающей.

Рассмотрим   функцию
$$
F(x)=\int\limits_1^x f(t)dt.
$$
Поскольку $F'(x)=f(x)>0$, то функция $F(x)$ возрастает вместе
с $x$. Если при $x\rightarrow +\infty$ функция $F(x)$ имеет конечный предел, то ряд \eqref{Ryad} сходится.
Если $F(x)$ возрастает до бесконечности вместе с $x$: $\lim\limits_{x\rightarrow +\infty}F(x)=+\infty$, то ряд \eqref{Ryad} расходится, но $F(x)$ позволяет судить о скорости роста частичной суммы
рассматриваемого ряда. Покажем это.

По формуле конечных приращений Лагранжа \eqref{Lag2} имеем
\begin{equation}\label{Rav}
	F(k+1)-F(k)=f(k+\theta),\qquad 0<\theta<1,
\end{equation}
а в силу строгого убывания функции $f(x)$ получим
\begin{equation}\label{Ner0}
	f(k)<F(k)-F(k-1)<f(k-1),
\end{equation}
и
\begin{equation}\label{Ner}
	f(k+1)<F(k+1)-F(k)<f(k).
\end{equation}
Отсюда получим
\begin{equation}\label{S0}
	0<f(k-1)-[F(k)-F(k-1)]<f(k-1)-f(k)
\end{equation}
и
\begin{equation}\label{S1}
	0<f(k)-[F(k+1)-F(k)]<f(k)-f(k+1).
\end{equation}
Просуммировав \eqref{S0} от $2$ до $n$ по $k$, а
\eqref{S1} от $1$ до $n$ по $k$, будем иметь
\begin{equation}\label{Ner00}
	0<	\sum\limits_{k=2}^n f(k-1)-[F(n)-F(1)]<f(1)-f(n)<f(1)
\end{equation}
и
\begin{equation}\label{Ner1}
	0<\sum\limits_{k=1}^n f(k)-[F(n+1)-F(1)]<f(1)-f(n+1)<f(1).
\end{equation}
Прибавив ко всем частям неравенства \eqref{Ner00} число $f(n)>0$, получим
\begin{equation}\label{Ner000}
	f(n)<\sum\limits_{k=1}^n f(k)-[F(n)-F(1)]<f(1).
\end{equation}
Рассмотрим последовательности
$$
\alpha_n=\sum\limits_{k=1}^n f(k)-F(n)
$$
и
$$
\beta_n=\sum\limits_{k=1}^n f(k)-F(n+1).
$$
Из \eqref{Ner000} и \eqref{Ner1} получим $f(n)-F(1)<\alpha_n$ и $\beta_n<f(1)-F(1)$,
то есть последовательность $\{\alpha_n\}$ ограничена снизу, а последовательность $\{\beta_n\}$  ограничена сверху.
Поскольку в силу \eqref{Ner}
$$
\alpha_{n+1}-\alpha_n=f(n+1)-[F(n+1)-F(n)]<0,
$$
$$
\beta_{n+1}-\beta_n=f(n+1)-[F(n+2)-F(n+1)]>0,
$$
то последовательность $\{\alpha_n\}$ убывает, а последовательность $\{\beta_n\}$ возрастает.
Значит, $\{\alpha_n\}$ сходится как убывающая и ограниченная снизу последовательность, а  $\{\beta_n\}$ сходится как возрастающая и ограниченная сверху последовательность.
Обозначим
$$
\lim\limits_{n\rightarrow+\infty}\alpha_n=A,
$$
$$
\lim\limits_{n\rightarrow+\infty}\beta_n=B.
$$
Тогда для любого $\varepsilon>0$ найдется $N\in\mathbb{N}$, что при $n>N$
$$
|\alpha_n-A|=\left|\sum\limits_{k=1}^n f(k)-F(n)-A\right|<\varepsilon,
$$
$$
|\beta_n-B|=\left|\sum\limits_{k=1}^n f(k)-F(n+1)-B\right|<\varepsilon
$$
или, полагая $S_n=\sum\limits_{k=1}^n f(k)$ и раскрывая модули,
$$
A+F(n)-\varepsilon<S_n<A+F(n)+\varepsilon,
$$
$$
B+F(n+1)-\varepsilon<S_n<B+F(n+1)+\varepsilon.
$$
Обозначая $C_1=A-\varepsilon$, $C_2=A+\varepsilon$, $C_3=B-\varepsilon$, $C_4=B+\varepsilon$, получим неравенства:
\begin{equation}\label{SL1}
	F(n)+C_1<S_n<F(n)+C_2,
\end{equation}
\begin{equation}\label{SL2}
	F(n+1)+C_3<S_n<F(n+1)+C_4.
\end{equation}
Учитывая, что $F(x)$ --- возрастающая функция, запишем и другие, более грубые, неравенства
\begin{equation}\label{SL3}
	F(n)+C_1<S_n<F(n+1)+C_2,
\end{equation}
\begin{equation}\label{SL4}
	F(n)+C_3<S_n<F(n+1)+C_4.
\end{equation}
Можно получить и варианты неравенств \eqref{SL3}--\eqref{SL4} с более аккуратными постоянными, которые также будут асимптотически точными, например,
$$
F(n)+C_5<S_n<F(n+1)+C_6.
$$
Это возможно в силу того, что величина $F(n+1)-F(n)$ является бесконечно малой, что следует из \eqref{Rav} и монотонного стремления $f(x)$ к нулю.

Описанным способом можно, например, получить неравенства \eqref{Eq23} для суммы гармонического ряда
$S_n=\sum\limits_{k=1}^n \dfrac{1}{k}$, что получается из  неравенств \eqref{SL2}
при выборе
$$f(x)=\frac{1}{x}, F(n+1)=\int\limits_1^{n+1}f(t)dt=\int\limits_1^{n+1}\frac{dt}{t}=\ln(n+1), C_3=0, C_4=1.$$

Далее мы будем употреблять название "неравенства Шлёмильха--Лемоннье"\  не только для сумм гармонического ряда, но и для других аналогичных рядов. Поясним, что будет иметься в виду под этим названием. Предположим, что из интегрального признака для рядов получены неравенства вида \eqref{SL1}--\eqref{SL2}. Пусть функции $G(x), H(x)$ эквивалентны на бесконечности интегральной функции $F(x)$, то есть выполнены предельные соотношения
$$
\lim_{x\to\infty}\frac{G(x)}{F(x)}=\lim_{x\to\infty}\frac{H(x)}{F(x)}=1.
$$
Тогда неравенствами Шлёмильха--Лемоннье мы будем называть любое из двусторонних неравенств для частичной суммы $S_n$ рассматриваемого ряда, которое имеет вид
\begin{equation}
G(n)+A_1 \leq S_n \leq H(n)+A_2,
\end{equation}
где $A_1, A_2$ --- некоторые постоянные.

Систематический вывод подобных асимптотических разложений и неравенств для рядов, интегралов и конечных сумм основан на использовании формулы суммирования Эйлера--Маклорена, см. \cite{Fikh2}.

\vskip 0.5 cm

\begin{pr}\label{extra}
	Для суммы
	$$
	R_n=\sum_{k=1}^n \frac{1}{\sqrt{k^2-1}}
	$$
выписать неравенства Шлёмильха--Лемоннье  вида \eqref{SL1}--\eqref{SL2}.
\end{pr}
{\it Решение.}
Имеем $f(x)=\frac{1}{\sqrt{x^2-1}}$,
$$
F(x)=\int\limits_1^x \frac{dt}{\sqrt{t^2-1}}=\ln(x+\sqrt{x^2-1}).
$$
Тогда неравенства примут вид
$$
\ln(n+\sqrt{n^2-1})+C_1<R_n<\ln(n+\sqrt{n^2-1})+C_2,
$$
$$
\ln(n+1+\sqrt{n(n+2)})+C_3<R_n<\ln(n+1+\sqrt{n(n+2)})+C_4.
$$

\vskip 0.5 cm

Теперь наметим план изучения еще одного класса сумм
$$
Q_n=\sum_{k=1}^n \frac{1}{k^\tau},\qquad 0<\tau<1.
$$
Соответствующий ряд называется обобщённым гармоническим рядом, при указанных значениях параметра $\tau$ он расходится.

\begin{zd}\hspace{-2mm}***\label{zd48}
 Изучите последовательности
\begin{equation}\label{Eq33}
a_n=Q_n-\frac{n^{1-\tau}}{1-\tau},
\end{equation}
\begin{equation}\label{Eq34}
b_n=Q_n-\frac{(n+1)^{1-\tau}}{1-\tau}.
\end{equation}
Докажите, что для $a_n, b_n$ существуют
$
\lim\limits_{n\rightarrow+\infty} a_n = C_1,
$
$
\lim\limits_{n\rightarrow+\infty} b_n = C_2
$
и $C_1=C_2$. Докажите неравенства типа Шлёмильха--Лемоннье для $Q_n$.
\end{zd}

\vskip 0.5 cm

В серьёзной математике выражение
\begin{equation}\label{Ri}
 \sum_{k=1}^\infty \frac{1}{k^s} = \zeta(s).
\end{equation}
называется "дзета--функцией Римана"\, от аргумента $s$. Мы показывали, что $\zeta(1)=+\infty$ (см. \eqref{GRRas}), этот факт выражает расходимость обычного гармонического ряда.
В соответствии с определением \eqref{Ri} указанный ряд сходится при условии ${\rm Re}(s)>1$.
Можно доказать, что
$$
\zeta(2)=\sum_{k=1}^\infty\frac{1}{k^2}=\frac{\pi^2}{6}\sim1,64493407,
$$
$$
\zeta(3)=\sum_{k=1}^\infty\frac{1}{k^3}\sim1,20205690.
$$
Значение $\zeta(2)$ пытались в своё время найти многие известные математики, это значение было найдено Эйлером, как и другие значения для чётных аргументов  $\zeta(2k)$. Точные значения $\zeta(3)$ и для всех остальных нечётных аргументов до сих пор не найдены в явном виде, более того неизвестно, являются ли они иррациональными или трансцендентными числами. Исключение составляет иррациональность числа  $\zeta(3)$, которая была доказана в 1978 г. Роже Апери, этот результат относится к числу наиболее значительных достижений в математике в двадцатом веке.

Можно также определить дзета--функцию Римана и для комплексных значений $0<{\rm Re}(s)<1$ при помощи так называемого аналитического продолжения. Самой известной из нерешенных до сих пор задач в математике является гипотеза Римана: все нули дзета--функции, определенной в полосе $0<{\rm Re}(s)<1$, лежат посередине этой полосы на прямой ${\rm Re}(s)=1/2$. Гипотеза Римана остается последней из великих задач, после того, как были решены проблема базисов (Пер Энфло, 1972), проблема Бибербаха (Луи де Бранж, 1984), доказана теорема Ферма (Эндрю Вайлс, 1993--1995).

Практически все остальные нерешенные задачи, связанные с натуральными числам и функциями от них, следуют из гипотезы Римана. Для них получены "условные"\, доказательства в предположении, что гипотеза Римана справедлива. Поэтому поиски доказательства этой гипотезы очень важны и ведутся до сих пор. Английский математик Гарольд Харди доказал, что на "критической"\, прямой ${\rm Re}(s)=1/2$ действительно лежат бесконечно много нулей $\zeta-$функции. На современных суперкомпьютерах проверено, что несколько десятков миллиардов первых нулей $\zeta-$функции, занумерованных в порядке возрастания их модулей, лежат на критической прямой ${\rm Re}(s)=1/2$. Поэтому справедливость гипотезы Римана не вызывает сомнений, хотя  мы и не умеем до сих пор ее строго доказывать.

Отметим, что идея использования дзета--функции и в целом методов теории функций комплексного переменного стали выдающимся вкладом Б.\,Римана в математику. Например, описание поведение дзета--функции на одной единственной прямой ${\rm Re} (s)=1$ позволило доказать основной закон распределения простых чисел. Сейчас дзета--функция активно используется при обосновании и оценках сложности компьютерных алгоритмов криптографии.

	\section{Доказательство неравенств методом математической индукции }

Математическая индукция --- мощный и изящный способ доказательства
некоторых типов неравенств: таких, которые
справедливы для всех натуральных чисел или для всех натуральных
чисел, начиная с некоторого числа.

\textbf{\textit{Метод математической индукции} } состоит в
следующем. Чтобы доказать, что некоторое утверждение верно для
любого натурального числа $n$, достаточно:

\begin{enumerate}
	\item   доказать это утверждение для $n=1$;
	
	\item   предположить его справедливость при $n=k$ и $k\ge 1$;
	
	\item   доказать, что оно верно при $n=k+1$.
\end{enumerate}

\textbf{\textit{Замечание.}} Методом математической индукции можно
доказывать утверждения, справедливые и при $n\ge m$, где $m>1$. В
ходе доказательства надо заменить \textit{первый шаг}: доказать
утверждение при $n=m$. Все остальное нужно оставить, как и прежде,
при необходимости пользуясь тем, что $n\ge m$.

\begin{pr}\label{MMIPr0}
	Доказать, что при любом натуральном $n$
	справедливо неравенство
	\begin{equation}\label{1}
		\left(\frac{n}{3}\right)^n<n!.
	\end{equation}
\end{pr}
{\it Решение.} Используем метод математической индукции.

1. Пусть $n=1$. Тогда верно, что
$$
\frac{1}{3}<1.
$$

2. Предположим, что верно неравенство для $n=k$:
\begin{equation}\label{2}
\left(\frac{k}{3}\right)^k<k!.
\end{equation}

3. Докажем неравенство для $n=k+1$:
$$
\left(\frac{k+1}{3}\right)^{k+1}<(k+1)!.
$$

Умножая обе части неравенства \eqref{2} на $\left(\frac{k+1}{3}\right)^{k+1}\,(k+1)$, получим
$$
\left(\frac{k}{3}\right)^k\,\left(\frac{k+1}{3}\right)^{k+1}\,(k+1)<(k+1)! \left(\frac{k+1}{3}\right)^{k+1}.
$$
После преобразования имеем
\begin{equation}\label{3}
\left(\frac{k+1}{3}\right)^{k+1}<(k+1)!\left[\frac{1}{3}\, \left(1+\frac{1}{k}\right)^{k}\right].
\end{equation}
Осталось доказать, что $\frac{1}{3}\, \left(1+\frac{1}{k}\right)^{k}<1$ или что $\left(1+\frac{1}{k}\right)^{k}<3$.
Сначала заметим, что при $k\geq 1$
$$
k!=k(k-1)...2\cdot 1\geq 2^{k-1}.
$$
Теперь, применяя формулу бинома Ньютона вида
$$
(x + y)^n =\sum\limits_{m=0}^nC_n^mx^{n-m}y^m,\qquad x,y\in\mathbb{R},\qquad n\in\mathbb{N},
$$
где $C_n^m=\frac{n!}{m!(n-m)!}$ -- биномиальный коэффициент,
получим
$$
\left(1+\frac{1}{k}\right)^{k}=1+\frac{1}{k}C_k^1+\frac{1}{k^2}C_k^2+...+\frac{1}{k^k}C_k^k=2+\sum\limits_{j=2}^k\frac{1}{k^j}C_k^j=
$$
$$
=2+\sum\limits_{j=2}^k\frac{1}{j!}\frac{k}{k}\frac{k-1}{k}...\frac{k-j+1}{k}\leq2+\sum\limits_{j=2}^k\frac{1}{2^{j-1}}=3-\frac{1}{2^{k-1}}<3.
$$
Тогда из \eqref{3}, получим
$$
\left(\frac{k+1}{3}\right)^{k+1}<(k+1)!\left[\frac{1}{3}\, \left(1+\frac{1}{k}\right)^{k}\right]<(k+1)!
$$
и
$$
\left(\frac{k+1}{3}\right)^{k+1}<(k+1)!.
$$

\begin{pr}\label{MMIPr1}
Доказать, что при любом положительном $a$ и при любом натуральном $n$
справедливо неравенство
\begin{equation}\label{GrindEQ__5_}
\underbrace{\sqrt{a+\sqrt{a+...+\sqrt{a} } } }_{n} \le \frac{1+\sqrt{4a+1} }{2},\qquad\forall n\in\mathbb{N},\qquad a>0.
\end{equation}
\end{pr}
{\it Решение.} Используем метод математической индукции.

1. Пусть $n=1$. Покажем, что
\[\sqrt{a} \le \frac{1+\sqrt{4a+1} }{2} \]
или
\[2\sqrt{a} -1\le \sqrt{4a+1} .\]
Это неравенство справедливо, поскольку
\[(2\sqrt{a} -1)^{2} =4a-4\sqrt{a} +1\le 4a+1.\]

2. Допустим, что
\begin{equation} \label{GrindEQ__6_} \underbrace{\sqrt{a+\sqrt{a+...+\sqrt{a} } } }_{k} \le \frac{1+\sqrt{4a+1} }{2} .
\end{equation}

3. Докажем неравенство для $k+1$ корня
\begin{equation} \label{GrindEQ__7_} \underbrace{\sqrt{a+\sqrt{a+...+\sqrt{a} } } }_{k+1} \le \frac{1+\sqrt{4a+1} }{2} .
\end{equation}

Прибавим $a$ к обеим частям равенства \eqref{GrindEQ__6_}, получим
\[a+\sqrt{a+\sqrt{a+...+\sqrt{a} } } \le \frac{1+\sqrt{4a+1} }{2} +a=\frac{2a+1+\sqrt{4a+1} }{2} .\]
Очевидно, что
\[\frac{2a+1+\sqrt{4a+1} }{2} =\left(\frac{1+\sqrt{4a+1} }{2} \right)^{2} ,\]
следовательно, неравенство \eqref{GrindEQ__7_} доказано, а,
следовательно, доказано и неравенство \eqref{GrindEQ__5_}.

В следующих задачах нужно методом математической индукции доказать, что для любого натурального числа $n$ (или для любого натурального числа $n$, начиная с некоторого) справедливы
неравенства.

\begin{zd}\label{zdMMI1}
	\begin{enumerate}[a)]
\item $n!<\left(\frac{n}{2} \right)^{n} ;$
		
\item $1+\frac{n}{2} \le \left(\frac{3}{2}
\right)^{n};$

	\item $n^{n+1}>(n+1)^n,\qquad n\geq 3;$
	
	\item $\frac{4^{n} }{n+1} <\frac{(2n)!}{(n!)^{2}};$
	
	\item $n!<n\left(\frac{n}{e}\right)^n.$
\end{enumerate}
\end{zd}

\begin{zd}\label{zdMMI2}
		\begin{enumerate}[a)]
			\item $\sqrt{n} \le
			\sqrt[{n}]{n!}<n+1;$

			\item $\frac{1}{\sqrt{n+1}}<2(\sqrt{n+1} -\sqrt{n}
			)<\frac{1}{\sqrt{n} } ;$
			
		\item $\sqrt{n}<1+\frac{1}{\sqrt{2}}+\frac{1}{\sqrt{3}}+...+\frac{1}{\sqrt{n}}<2\sqrt{n};$
	
	\item  $\frac{n}{2}\leq 1+\frac{1}{2}+...+\frac{1}{2^n-1}<n.$
	
\end{enumerate}
\end{zd}

\begin{zd}\label{zdMMI3}
\begin{enumerate}[a)]
	\item  $1+\frac{1}{2^{2} } +\frac{1}{3^{2} }
+...+\frac{1}{n^{2} } \le 2-\frac{1}{n};$

\item $\left(1-\frac{1}{\sqrt{2} } \right)\left(1-\frac{1}{\sqrt{3} }
\right)...\left(1-\frac{1}{\sqrt{n} } \right)<\frac{2}{n^{2} }
,\quad \quad n>1;$

\item $\left({\rm \; 1}+{\rm \; }\frac{{\rm 1}}{{\rm
		1}^{{\rm 3}} } \right)\left({\rm 1\; }+{\rm \; }\frac{{\rm
		1}}{{\rm 2}^{{\rm 3}} } \right)\left({\rm 1\; }+{\rm \;
}\frac{{\rm 1}}{{\rm 3}^{{\rm 3}} } \right){\rm ...}\left({\rm 1\;
}+{\rm \; }\frac{{\rm 1}}{{\rm n}^{{\rm 3}} } \right){\rm \;
}<{\rm \; 3\; -}\frac{{\rm 1}}{{\rm n}}.$
\end{enumerate}
\end{zd}




\begin{zd}\label{zdMMI4}
\begin{enumerate}[a)]
	\item  $2!\cdot 4!\cdot 6!\cdot...(2n)!\geq
((n+1)!)^n;$
\item  $ \frac{1}{2} \cdot \frac{3}{4} \cdot \frac{5}{6}
...\frac{2n-1}{2n} \le \frac{1}{\sqrt{3n+1} } ;$
\item  $\frac{1\cdot 3\cdot 5\cdot ...\cdot (2n-1)}{n!}
<\frac{2^{n} }{\sqrt{2n+1} }.$
\end{enumerate}
\end{zd}

\begin{zd}\label{zdMMI5} Пусть $a,b\in {\rm \mathbb{R}}$, такие что $ a+b>0$ и $a\ne b$. Показать, что при $n>1$ справедливо неравенство
$$(a+b)^{n} <2^{n-1} (a^{n} +b^{n} ).$$
\end{zd}

\begin{zd}\label{zdMMI6} Пусть $x_{1} ,...,x_{n} \in
	{\mathbb{R}}$ и $x_{1} >0,...,$ $x_{n} >0$. Докажите неравенства.
	\begin{enumerate}[a)]
\item  $(\sqrt{x_{1} } +...+\sqrt{x_{n} } )^{2} \le
n(x_{1} +...+x_{n} )$;

\item  $\frac{n^{3} }{(x_{1} +...+x_{n} )^{2} } \le
\frac{1}{x_{1}^{2} } +...+\frac{1}{x_{n}^{2} } $;

\item  $nx_{1} ...x_{n} \le x_{1}^{n} +...+x_{n}^{n} $;

\item  $n\le x_{1} +...+x_{n} ,$ при дополнительном ограничении  $x_{1} ...x_{n}
=1$;

\item  $\frac{1}{\frac{1}{x_{1} } +...+\frac{1}{x_{n} }
} \le \sqrt[{n}]{x_{1} ...x_{n} } \le \frac{x_{1}
	+...+x_{n} }{n} $.
\end{enumerate}
\end{zd}

Последние неравенства являются частными случаями неравенств между различными степенными средними, см. далее.

\section{Неравенства для числа $e$ и связанных с ним функций}

Теперь перейдем еще к одному циклу задач. Он связан с замечательным числом "$e$"\  в математике, основанием натуральных логарифмов. Напомним, что это число определяется как предел последовательности
$$
\lim_{n\to \infty} \left(1+\frac{1}{n}\right)^n=e,
$$
кроме того, справедливы равенства
$$
e=\lim_{n\to \infty} \left(1+\frac{1}{n}\right)^{n+1}=\lim_{x\to {+\infty}} \left(1+\frac{1}{x}\right)^x=
\lim_{x\to {-\infty}} \left(1+\frac{1}{x}\right)^x.
$$

Рассмотрим подробнее несколько задач на определение числа $e$.

\begin{zd}
Докажите, что последовательность
$$
a_n=\left(1+\frac{1}{n}\right)^n
$$
монотонно возрастает и ограничена сверху.
\end{zd}

{\it Указание:} для доказательства монотонности можно использовать один из трёх вспомогательных результатов: бином Ньютона, неравенство Бернулли или неравенство для средних.

\begin{zd}
Докажите, что последовательность
$$
b_n=\left(1+\frac{1}{n}\right)^{n+1}
$$
монотонно убывает и ограничена снизу.
\end{zd}

Таким образом, по теореме Вейерштрасса существуют два предела
$$
e_1=\lim_{n\to \infty} a_n, \ \ \  e_2=\lim_{n\to \infty} b_n.
$$

\begin{zd}
Докажите, что выполняется равенство $e_1=e_2$.
\end{zd}
Таким образом, можно обозначить общую величину двух пределов через $e=e_1=e_2$. Это и есть знаменитое число "$e$".

\begin{zd}
Докажите, что выполнены неравенства $2,25 < e < 4 $.
\end{zd}

Современное обозначение для числа е было введено Л.\,Эйлером, а его приближённое значение впервые вычислил Я.\,Бернулли.

На самом деле, приближённое значение числа $e$ такое:
$$
e\approx 2,718281828459045...
$$
Для запоминания можно обратить внимание на дважды повторяемый набор цифр 1828 --- это год рождения Льва Толстого (который называл себя Лёв) или год написания Бетховеном Лунной сонаты (которую Бетховен никак не назвал), но математики обычно наоборот вспоминают эти даты из числа $e$. Число $e$ является не только иррациональным, но даже трансцендентным. Это означает, что оно не является корнем никакого многочлена с целыми коэффициентами. Иррациональность числа е установить нетрудно, см. \cite{Fikh}, трансцендентность была доказана Ш.\, Эрмитом.

\begin{zd}\label{zd50}
Докажите, что при $x>0$
\begin{equation}\label{Eq35}
\bigg(1+\frac{1}{x}\bigg)^x<e<\bigg(1+\frac{1}{x}\bigg)^{x+1}.
\end{equation}
\end{zd}

Оказывается, правую часть неравенства \eqref{Eq35} можно уточнить.
Для изучения этого неравенства введем функцию $\varepsilon(x)$ из соотношения
$$
\bigg(1+\frac{1}{x}\bigg)^{x+\varepsilon(x)}=e,\qquad x>0,\quad \text{следовательно},
$$
$$
\varepsilon(x)=\frac{1}{\ln\left( 1+\frac{1}{x}\right)}-x,\qquad x>0.
$$

Функция $\varepsilon(x)$ определяет отклонение величины $(1+\frac{1}{x})^{\frac{1}{x}}$ от числа $e$.

\begin{zd}\label{zd52} Докажите, что при $x>0$
$$
0<\varepsilon(x)<1.
$$
\end{zd}

\begin{zd}\label{zd53}
 Докажите, что $\varepsilon(x)$ монотонно возрастает на полуоси $(0,\infty)$, причем
$$
0<\varepsilon(x)<\frac{1}{2}
$$
\end{zd}

{\it Указание:} покажите, что $\varepsilon'(x)>0$ и
$$
\lim_{x\rightarrow\infty}\varepsilon(x)=\frac{1}{2},
$$
и значит, применима теорема \ref{teo2}. Для вычисления предела дважды примените правило Бернулли--Лопиталя.

\begin{zd}\label{zd54}
Докажите, что при $x>0$
\begin{equation}\label{Eq37}
\text{a)}\ \  e<\bigg(1+\frac{1}{x}\bigg)^{x+\frac{1}{2}}.
\end{equation}
\begin{equation}\label{Eq38}
\text{b)}\ \  \frac{1}{x+\frac{1}{2}}<\ln\left( 1+\frac{1}{x}\right).
\end{equation}
\end{zd}

Из Задачи \ref{zd54} следует, что постоянную $\frac{1}{2}$ в \eqref{Eq37}--\eqref{Eq38} нельзя заменить никаким меньшим числом. Следовательно постоянная $\frac{1}{2}$ --- точная в этих неравенствах.

\begin{zd}\label{zd55}
 Пусть $x\geqslant1$. Тогда
$$
\varepsilon_{2n-1}(x)\leqslant\varepsilon(x)\leqslant\varepsilon_{2n}(x),\quad n=1,2,\ldots, \text{где}
$$
$$
\varepsilon_n(x)=\frac{1}{\sum\limits_{j=1}^n\frac{(-1)^{j-1}}{jx}}-x.
$$
\end{zd}

{\it Указание:} используйте формулу Тейлора для логарифма и теорему Лейбница о знакопеременном ряде.

\begin{zd}\label{zd56}
Докажите, что
$$
\text{a)}\ \  \, e<\bigg(1+\frac{1}{x}\bigg)^{x+\varepsilon_{2n}(x)},\quad x\geqslant1;
$$
$$
\text{b)}\ \  \, e>\bigg(1+\frac{1}{x}\bigg)^{x+\varepsilon_{2n-1}(a)},\quad x\geqslant a.
$$
\end{zd}

Теперь используем идею из Примера \ref{pr6}.

\begin{zd}\label{zd57} Докажите, что
$$
e>\bigg(1+\frac{1}{x}\bigg)^{x+c_k},\quad x\geqslant k, k=1,2,3,
$$
с постоянными $c_1=\frac{1}{5}, \, c_2=\frac{2}{5}, \, c_3=\frac{21}{47}.$
\end{zd}

\begin{zd}\label{zd58}
Докажите, что при $x\geqslant1$
$$
\text{a)}\ \  \, e>\left(1+\frac{1}{x}\right)^{x+\frac{1}{2}-\frac{x+2}{2(6x^2-9x+2)}}
$$
$$
\text{b)}\ \  \, e>\left(1+\frac{1}{x}\right)^{x+\frac{1}{2}-\frac{2x^2-2x-9}{2(12x^3-6x^2+4x-9)}}
$$
\end{zd}

\begin{zd}\hspace{-2mm}**\label{zd59}
Найдите как можно более точно такое число $a$, что $\forall x\geqslant a$
$$
e>\bigg(1+\frac{1}{x}\bigg)^{x+\frac{2}{5}}.
$$
\end{zd}

{\it Указание:} $a \approx 0,413053$, это корень уравнения
$$
\left(1+\frac{1}{x}\right)^{x+0,4}=e,\ \  \text{или}\ \  \varepsilon(x)=0,4.
$$

Рассмотрим уравнение из последней задачи
$$
f(x)=\frac{1}{\ln(1+x)}-0,4=x.
$$
Метод итераций для такого уравнения вблизи $a$ расходится, условие $|f'(a)|{<}1$ не выполнено. Известно, что в этом случае сходится итерационный процесс для обратной функции $f^{-1}$:
$$
x_{n+1}=\frac{1}{{\rm exp}\left(\frac{1}{x_n+0,4}\right)-1},
$$
см. \cite{VM}.
Однако, вычисления показывают, что для достижения указанной точности $a$ при $x_0=1$ требуется около ста итераций. Практически вычисления вообще невозможно довести до конца из-за накопления погрешности.

С другой стороны, используя для ускорения сходимости "$\lambda$--метод"\, cм.  \cite{VM}, получаем процесс
$$
x_{n+1}=x_n+\lambda\left( \frac{1}{\ln\big(1+\frac{1}{x_n}\big)}-x_n-0,4\right)
$$
при $\lambda=-7,47, \, x_0=1$, который вычисляет приведенное выше значение $a$ за 4--5 итераций. Этот пример демонстрирует силу "$\lambda$--метода"\, ускорения сходимости.

 {\it Вопрос:} почему мы выбрали именно такое значение $\lambda$?

Аналогично исследуется случай $x<0$.

\begin{zd}\label{zd60}
Докажите неравенства
$$
\text{a)}\quad \bigg(1+\frac{1}{x}\bigg)^{x+\frac{1}{2}}<e<\bigg(1+\frac{1}{x}\bigg)^{x+1},\quad x<-1,
$$
$$
\text{b)}\quad e<\bigg(1+\frac{1}{x}\bigg)^{x+0,6},\quad x<-1,413053
$$
\end{zd}
{\it Указание:} используйте равенство $\varepsilon(-x)=1-\varepsilon(x-1),$ $x>1$.

\begin{zd}\label{zd61}
Докажите, что при $x\geqslant1$
$$
e^{\frac{1}{2x}-\frac{1}{3x^2}}<\frac{e}{(1+\frac{1}{x})^x}<e^{\frac{1}{2x}}.
$$
\end{zd}

 Задача \ref{zd61} --- это существенное уточнение задачи \ref{zd50}.

\begin{zd}\label{zd62}
Исследуйте сходимость ряда
$$
1+\ln2+\ln2\cdot\ln\left( 1+\frac{1}{2}\right)^2+\ln2\cdot\ln\left( 1+\frac{1}{2}\right)^2\cdot\ln\left( 1+\frac{1}{2}\right)^3+\ldots
$$
\end{zd}
{\it Указание:} признаки Коши и Даламбера не дают ответа.
Примените признак Раабе.

Очевидно, что экспоненциальная функция является положительной на всей действительной оси. Оказывается, что для очень похожей на неё функции положительность является нерешённой проблемой.

\begin{zd}\hspace{-2mm}***  Докажите, что следующая функция является положительной на всей действительной оси
$$
f(x) = \sum_{k=0}^{\infty} \frac{x^k}{\sqrt{n!}}>0,\ \  x\in \mathbb{R}.
$$
\end{zd}

Разумеется, последнее неравенство является содержательным только для отрицательных $x$.

\begin{zd}\label{zd63}
 Докажите неравенства для логарифмической функции
\begin{enumerate}[a)]
\item $\displaystyle\ln^2\Big(1+\frac{1}{x}\Big)<\frac{1}{x(x+1)}$, \quad $x>0$,
\item $\displaystyle\ln\Big(1+\frac{1}{x}\Big)^x\cdot\ln\Big(1+\frac{1}{x}\Big)^{x+1}<1$, \quad $x>0$,
\item $\displaystyle\frac{\ln\Big(1+\displaystyle\frac{1}{y}\Big)-\ln\Big(1+\frac{1}{x}\Big)}{x-y}>\ln\Big(1+\frac{1}{x}\Big)\cdot\ln\Big(1+\frac{1}{y}\Big)$, \quad $x>y$,
\item $\displaystyle\ln^2\Big(1+\frac{1}{x}\Big)\leq\frac{1}{x(x+1)}-\frac{1}{x^2}\cdot\sum\limits_{k=0}^{2n+1}\frac{(-1)^k\cdot(1-A_k)}{x^k}$, \quad \mbox{$\text{где}\quad x\geq 1$,}
\end{enumerate}
$\displaystyle A_k=\underset{\substack{i+j=k\\ i\geq0,j\geq0}}{\sum}\frac{1}{(i+1)\cdot(j+1)}$,\\
и, если заменить $2n+1$ на $2n$, то будет справедливо обратное неравенство.
\end{zd}
{\it Указание:} рассмотрите функцию $\displaystyle\varepsilon(x)=\frac{1}{x(x+1)}-\ln^2\Big(1+\frac1x\Big)$ и разложите ее в ряд по степеням $\displaystyle\frac1n$.

\begin{zd}\label{zd631}
	Докажите неравенства для логарифмических функций
\begin{enumerate}
\item[a)] $\displaystyle\displaystyle\ln^2\Big(1+\frac{1}{x}\Big)<\frac{1}{x(x+1)}<\ln^{1,99}\Big(1+\frac{1}{x}\Big)$, \quad $x\geq3$,
\item[b)] $\displaystyle\ln^2\Big(1+\frac{1}{x}\Big)<\frac{1}{x(x+1)}<\ln^{1,99}\Big(1+\frac{1}{x}\Big)$, \quad $x\geq10$.
\end{enumerate}
\end{zd}

Теперь мы приведем ряд открытых проблем, решения которых авторам неизвестны. Это задачи \ref{zd64} и \ref{zd65}. Проверьте свои силы, доказав эти неравенства аналитически или на компьютере.

\begin{zd}\hspace{-2mm}***\label{zd64}
Докажите, что при $x>0$ выполнены неравенства для логарифмических функций
$$
\displaystyle\ln^3\Big(1+\frac1x\Big)<\frac1{x(x+1)(x+\frac12)},
$$
$$
\displaystyle\ln^4\Big(1+\frac1x\Big)<\frac1{x(x+1)(x+\frac12)(x+0,35)},
$$
$$
\displaystyle\ln^n\Big(1+\frac1x\Big)<\frac1{\Big(x+\displaystyle\frac n2\Big)x^{n-1}}, \qquad n\geq 2.
$$
\end{zd}

Отметим, что задача о разложении в ряд функции \mbox{$\ln^n(1+y)$} приводит к числам Стирлинга, находящим многочисленные применения.

\begin{zd}\hspace{-2mm}***\label{zd65}
Определим функцию $\varepsilon(z)$ при комплексных значениях $z\in \mathbb{C}[-1,0]$ по формулам
$$
\Bigg(1+\frac1z\Bigg)^{z+\varepsilon(z)}=e, \qquad \varepsilon(z)=\frac1{\ln\Big(1+\displaystyle\frac1z\Big)}-z.
$$
Верно ли, что $\displaystyle|\varepsilon(z)|\leq\frac12$?
\end{zd}

Ранее мы доказали, что это верно, если $z=x>0$, $x\in \mathbb{R}$.

Приведем несколько задач повышенной трудности, связанных с гармоническим рядом и дзета--функцией Римана \eqref{Ri}.

\begin{zd}\hspace{-2mm}*\label{zd66}
Докажите, что
$$
\displaystyle\lim\limits_{n\rightarrow+\infty}\sum\limits_{k=1}^n\frac{(-1)^{k+1}}{k}=\ln 2.
$$
\end{zd}
Таким образом, если в гармоническом ряде поменять знак у каждого второго члена, то он уже будет сходиться.

\begin{zd}\hspace{-2mm}**\label{zd67}
 В гармоническом ряде $S_n$ отбросили все члены вида $\displaystyle\frac1n$, если в десятичной записи числа $n$ встречается цифра 7.
Докажите, что ряд из оставшихся членов будет сходящимся. Оцените сверху его сумму.
\end{zd}

\begin{zd}\hspace{-2mm}**\label{zd68}
 Докажите равенство
$$
\displaystyle S_n=c+\ln(n)+\frac1{2n}-\sum\limits_{k=2}^\infty\frac{A_k}{n(n+1)\ldots(n+k-1)},
$$
где числа $A_k$ вычисляются по формуле
$$
\displaystyle A_k=\frac1k\int\limits_0^1x(1-x)(2-x)(3-x)\ldots(k-1-x)dx.
$$

В частности,
$$
\displaystyle A_2=\frac1{12}, \quad \displaystyle A_3=\frac1{12}, \quad \displaystyle A_4=\frac{19}{80}, \quad \displaystyle A_2=\frac9{20}.
$$
\end{zd}

 Задача \ref{zd68} --- это источник бесконечного улучшения оценок и предельных соотношений для $S_n$.

\begin{zd}\label{zd69}
Вычислите пределы
$$
\lim\limits_{n\rightarrow+\infty}n(S_n-C-\ln(n)),
$$
$$
\displaystyle\lim\limits_{n\rightarrow+\infty}n^2\left( S_n-C-\ln(n)-\frac1{2n}\right).
$$
\end{zd}

\begin{zd}\hspace{-2mm}*\label{zd70}
Выведите из равенства
$$
\displaystyle\sum\limits_{k=1}^\infty\frac1{k^a}=\zeta(a),
$$
при ${\rm Re}(a)>1$ равенство
$$
\displaystyle\sum\limits_{k=1}^\infty\frac{(-1)^{k+1}}{k^a}=(1-2^{1-a})\cdot\zeta(a).
$$
\end{zd}

\begin{zd}\hspace{-2mm}**\label{zd71}
Числа Бернулли $B_{2n}$ определяется как коэффициенты разложения в ряд Тейлора функции:
$$
\displaystyle\frac{x}{e^x-1}=1-\frac x2+\sum\limits_{k=1}^\infty\frac{B_{2k}\cdot x^{2k}}{(2k)!}, \qquad |x|<2\pi.
$$

Докажите, что
$$
\displaystyle\sum\limits_{k=1}^\infty\frac1{k^{2n}}=\frac{2^{2n-1}\cdot\pi^{2n}}{(2n)!}\cdot|B_{2n}|,
$$
$$
\displaystyle\sum\limits_{k=1}^\infty\frac{(-1)^{k+1}}{k^{2n}}=\frac{(2^{2n-1}-1)\cdot\pi^{2n}}{(2n)!}\cdot|B_{2n}|.
$$
\end{zd}

Эти формулы и были выведены Л.\,Эйлером.

\begin{zd}\hspace{-2mm}*\label{zd72}
Докажите равенства:
$$
\displaystyle\sum\limits_{k=1}^\infty\frac1{k^2}=\frac{\pi^2}6, \qquad\qquad \sum\limits_{k=1}^\infty\frac1{(2k-1)^2}=\frac{\pi^2}8,
$$
$$
\displaystyle\sum\limits_{k=1}^\infty\frac{(-1)^{k+1}}{k^2}=\frac{\pi^2}{12}, \qquad\qquad \sum\limits_{k=1}^\infty\frac1{(2k-1)^4}=\frac{\pi^4}{96}.
$$
\end{zd}
{\it Указание:} используйте равенство Парсеваля для рядов Фурье простейших функций.

В заключение цикла задач на исследование конкретных функций рассмотрим еще один метод получения логарифмических неравенств, основанный на разложении логарифма в непрерывную дробь.
Так называется представление вида
$$
        \ln(1+x)=\frac{x}{1+\displaystyle\frac{1\cdot x}{2+\displaystyle\frac{1\cdot x}{3+\displaystyle\frac{2^2\cdot x}{4+\displaystyle\frac{2^2\cdot x}{5+\displaystyle\frac{3^2\cdot x}{6+\displaystyle\frac{3^2\cdot x}{7+\ldots}}}}}}}.
$$
Это разложение в бесконечную (или "цепную"\,) дробь было известно еще Эйлеру и Гауссу. Термин "непрерывная дробь"\, ввел Валлис.

Если оборвать непрерывную дробь, отбросив ее бесконечную часть $n$--го знака плюс, то получится $n$--ая подходящая дробь вида
$$
\displaystyle R_n(x)=\frac{P_n(x)}{Q_n(x)},
$$
где $P_n(x)$ и $Q_n(x)$ --- некоторые многочлены. Можно показать, что при всех $x>-1$ справедливо равенство
$$
\lim\limits_{n\rightarrow+\infty}R_n(x)=\ln(1+x),
$$
то есть непрерывная дробь сходится к логарифму при $x>-1$. В этом проявляется большое преимущество непрерывных дробей перед рядами Тейлора.
Напомним, что ряд Тейлора функции $\ln(1+x)$ сходится лишь при $-1<x\leq1$.

\begin{zd}\label{zd73}
 Вычислите первые подходящие дроби
$$
R_1=x, \quad R_2=\frac{2x}{x+2}, \quad R_3=\frac{x(x+6)}{4x+6},
$$
$$
R_4=\frac{x(3x+6)}{x^2+6x+6}, \quad R_5={x(x^2-13x+30)}{9x^2+28x+30}.
$$
\end{zd}

\begin{zd}\hspace{-2mm}**\label{zd74}
 Докажите, что для любого $k=0,1,2,\ldots$
\begin{equation}\label{Eq39}
R_{2k}(x)<\ln(1+x)<R_{2k+1}(x), \quad x>0.
\end{equation}
\end{zd}

{\it Указание:} обозначить отброшенную часть непрерывной дроби через $\varepsilon>0$.

Можно также доказать (попробуйте!!!), что каждое последующее неравенство \eqref{Eq39} точнее предыдущего,
то есть для любого $x>0$ и для любого $k\in \mathbb{N}$
$$
R_{2k+3}(x)<R_{2k+1}(x), \quad R_{2k}(x)<R_{2k+2}(x).
$$

Таким образом, используя метод непрерывных дробей, мы получаем бесконечную последовательность все более и более точных двусторонних неравенств.
При этом уже на первом шаге получаются неравенства \eqref{Eq05} и \eqref{Eq10}. Делая замену $x\rightarrow\frac1x$,  получим из равенства для $\displaystyle\ln\Big(1+\frac1x\Big)$.

\begin{zd}\hspace{-2mm}**\label{zd75}
Докажите, что для подходящей дроби $R_n$ степень многочлена $P_n(x)$ равна $\displaystyle\left[\frac{n+1}2\right]$, а степень многочлена $Q_n(x)$ равна $\displaystyle\left[\frac n2\right]$, где $[t]$ -- есть целая числа $t$.
\end{zd}

\begin{zd}\hspace{-2mm}***\label{zd76}
	Найдите формулы для многочленов $P_n(x)$, $Q_n(x)$.
\end{zd}

Решение последней задачи авторам неизвестно.

Теперь выясним общую закономерность, которая управляет задачами типа \ref{zd20}--\ref{zd25}, \ref{zd28}--\ref{zd36}, \ref{zd48}.

В последующих задачах мы продолжаем использовать обозначения:

$f(x)$ --- функция, определенная на $[1,\infty)$, непрерывно дифференцируемая, монотонно убывающая, неотрицательная;

$F(x)=\int\limits_1^xf(t)dt$  --- её первообразная;

$S(n)=\sum\limits_{k=1}^n\,f(k)$.

\begin{zd}\label{zd77}
Доказать, что $F(x)$ монотонно возрастает.
\end{zd}

\begin{zd}\label{zd78}
Докажите справедливость двустороннего неравенства
\begin{equation}\label{Eq40}
F(n)\leq S(n)\leq F(n)+a_1-a_{n+1}.
\end{equation}
\end{zd}

\begin{zd}\label{zd79}
Покажите, что это дает для сумм гармонического ряда $S_n$.
\end{zd}

\begin{zd}\label{zd80}
Что это дает для сумм $P_n$?
\end{zd}

\begin{zd}\label{zd81}
Докажите, что ряд
$$
\sum\limits_{k=1}^\infty f(k)
$$
сходится тогда и только тогда, когда конечен предел
$$
\lim\limits_{n\rightarrow\infty}F(n).
$$
\end{zd}

\begin{zd}\label{zd82}
Докажите, что последовательность
$$
a_n=S(n)-F(n)
$$
монотонно убывает.
\end{zd}

\begin{zd}\label{zd83}
 Докажите, что последовательность
$$
b_n=S(n)-F(n+1)
$$
монотонно возрастает.
\end{zd}

\begin{zd}\label{zd84}
Докажите, что существует предел $\lim a_n=C_1$.
\end{zd}

\begin{zd}\label{zd85}
 Докажите, что существует предел $\lim b_n=C_2$.
\end{zd}

\begin{zd}\label{zd86}
Докажите, что если
$$
\lim\limits_{n\rightarrow+\infty}f(n)=0, \qquad C_1=C_2=C.
$$
\end{zd}

Эту величину $C$ назовем постоянной Эйлера данного ряда или данной последовательности.

\begin{zd}\label{zd87}
Докажите, что
\begin{equation}\label{Eq41}
0\leq C\leq a_1.
\end{equation}
\end{zd}

\begin{zd}\label{zd88}
Приведите пример функции $f(x)$, удовлетворяющей необходимым условиям, для которой $C_1\neq C_2$.
\end{zd}

\begin{zd}\hspace{-2mm}***\label{zd89}
Может ли в \eqref{Eq41} достигаться равенство?
\end{zd}

Авторы не знают ответа к задаче \ref{zd89}.

\begin{zd}\label{zd90}
Изучите сумму
$$
\displaystyle\sum\limits_{k=1}^n\frac1{k\ln k}.
$$
\end{zd}

\begin{zd}\hspace{-2mm}*\label{zd91}
Изучите самостоятельно случай возрастающей $f(x)$.
\end{zd}

\begin{zd}\label{zd92}
Изучите сумму
$$
\sum\limits_{k=1}^nk^a.
$$
\end{zd}

\begin{zd}\hspace{-2mm}*\label{zd93}
Изучите сумму
$$
\sum\limits_{k=1}^ne^{k^2}.
$$
\end{zd}

\begin{zd}\label{zd94}
Докажите неравенство
$$
\int\limits_0^xe^{t^2}dt>\frac1{2x}e^{x^2}, \qquad x>1.
$$
\end{zd}

Покажите, что при $x\rightarrow{+\infty}$ две части последнего неравенства эквивалентны, а ограничение $x>1$ можно улучшить.

Основной вывод из задач этого цикла: если функция $f(x)$ удовлетворяет перечисленным выше требованиям, то
$$
\sum\limits_{k=1}^nf(k)=\int\limits_1^nf(x)dx+C+L_n,
$$
где $C=C(f)$ соответствующая константа Эйлера, $L_n-$бесконечно малая. Кроме того, функция удовлетворяет неравенствам типа Шлёмильха--Лемоннье.

В заключение сформулируем несколько проблем, представляющих интерес и для специалистов.

\begin{zd}\hspace{-2mm}**\label{zd95}
Связана ли введенная нами постоянная Эйлера $C$ с постоянной в формуле Эйлера--Маклорена?
\end{zd}

\begin{zd}\hspace{-2mm}***\label{zd96}
Известно, что для простых чисел $P_n$
$$
\sum\limits_{k=1}^n\frac1{P_n}\approx\ln\ln n.
$$
Являются ли последовательности
$$
\sum\limits_{k=1}^n\frac1{P_n}-\ln\ln n, \quad \sum\limits_{k=1}^n\frac1{P_n}-\ln\ln(n+1)
$$
монотонными и выполняются ли для них неравенства типа Шлёмильха--Лемоннье?
\end{zd}

Оценки в этой задаче ведут к оценкам для $P_n$.

Рассмотрим дзета--функцию для значений $a>1$. Тогда очевидно, что
$$
\lim\limits_{n\rightarrow+\infty}\left(\sum\limits_{k=1}^n\frac1{k^a}-\frac1{a-1}\left(1-\frac1{n^{a-1}}\right)\right)=\zeta(a)-\frac1{a-1}.
$$

\begin{zd}\label{zd97}
Пусть теперь $0<a<1$. Будет ли формула
$$
\lim\limits_{n\rightarrow+\infty}\left(\sum\limits_{k=1}^n\frac1{k^a}-\frac{n^{1-a}}{1-a}\right)
$$
определять аналитическое продолжение $\zeta(a)$? Возможно ли аналогичное продолжение для $a<0$? Можно ли вычислить константы Эйлера $C(a)$ и справедливы ли неравенства типа Шлёмильха--Лемоннье?
\end{zd}

\begin{zd}\hspace{-2mm}***\label{zd98}
 Для нулей функции Бесселя $J_v(x)$ справедлива асимптотическая формула Макмагона
$$
\gamma_{v,n}=\left(n+\frac12v-\frac14\right)\pi+o(1).
$$

Является ли последовательность
$$
a_n=\gamma_{v,n}-\left(n+\frac12v-\frac14\right)\pi
$$
монотонной? Есть ли неравенства типа Шлёмильха--Лемоннье?
\end{zd}

\begin{zd}\hspace{-2mm}***\label{zd99}
 Задача Штурма--Лиувилля для уравнения
$$
y''(x)-g(x)y(x)=\lambda y(x), \quad a<x<b,
$$
при подходящих двухточечных краевых условиях порождает спектр с асимптотикой
$$
\lambda_n=n+o(1).
$$

Является ли последовательность $a_n=\lambda_n-n$ монотонной? Есть ли неравенства типа Шлёмильха--Лемоннье?
\end{zd}

Те же вопросы можно поставить и для спектральной функции, они более подходят к рассмотренному методу.

 Задача \ref{zd69} порождает новый ряд вопросов, когда величина $S(n)-F(n)-C$ раскладывается в разумный ряд по степеням $1/n$.

\begin{zd}\hspace{-2mm}***\label{zd100}
\begin{enumerate}[a)]
	\item Являются ли последовательности задачи \ref{zd69} монотонными? Какие условия на $f(n)$ гарантируют это?
\item Можно ли оценить эти последовательности? Такие весовые оценки было бы естественно назвать неравенствами Шлёмильха--Лемоннье второго рода, третьего рода (ранга) и т. д.
\end{enumerate}
\end{zd}

\section{Средние значения и неравенства для них}

\subsection{Различные виды средних значений}\label{sre}

Средние величины первоначально возникли в Древней Греции при решении задач на пропорции в терминах отношения отрезков и других геометрических задач, см. \cite{Sit2}. Так были введены простейшие средние: арифметическое, геометрическое, квадратичное и гармоническое
$$
A(x,y)=\frac{x+y}{2},\qquad G(x,y)=\sqrt{x,y},\qquad Q(x,y)=\sqrt{\frac{x^2+y^2}{2}},\quad x,y \ge 0,
$$
$$
H(x,y)=A\left(\frac{1}{x},\frac{1}{y}\right)=\frac{2xy}{x+y}, \qquad x,y>0.
$$

Наблюдения за приведёнными выше простейшими средними приводят к выводу, что все они обладают некоторыми характерными общими свойствами: однородны, симметричны, монотонны, а также  при $x=y$ равны $x$. Это приводит к идее определять абстрактные средние при помощи соответствующих аксиом.

Далее материал излагается в основном на основе \cite{Sit1, Sit2}, а также работ \cite{Sit3}--\cite{Sit15}.

\textbf{ Определение.}
\textit{Абстрактным средним} двух неотрицательных чисел   называется
число $M(x,y)$, удовлетворяющее следующим аксиомам:

1) (свойство несмещённости)
$$
M(x,x) = x ,
$$

2) (свойство однородности)
$$
M(\lambda x, \lambda y) = \lambda \ M(x,y), \lambda > 0,
$$

3) (свойство монотонности по обоим аргументам)
$$
x_{2} > x_{1} \Rightarrow M(x_{2}, y) > M(x_{1}, y),\
 y_{2} > y_{1}\Rightarrow M(x, y_{2}) > M(x, y_{1}),
$$

4)  (свойство симметричности)
$$
M(x,y) = M(y,x) .
$$

Очевидно, что если выполнено свойство симметричности, то тогда
достаточно монотонности лишь по одному аргументу. Из приведённых
выше  вытекает также важное свойство промежуточности
$$
\min(x,y)\le M(x,y) \le \max (x,y),
$$
которое часто включают в число основных (аналогично тому, как
включают в число аксиом нормы неотрицательность, хотя она следует из
других аксиом). Обычно подразумевается также непрерывность по обоим
аргументам, иногда для получения более детальных результатов нужны
дополнительно некоторая гладкость (например, если неравенства между
средними доказываются с использованием одной или нескольких
производных) и определённая модификация свойства аналитичности.

Далее нам также потребуется  понятие сопряженного среднего.

\textbf{ Определение.}  Сопряженным к абстрактному среднему
называется величина
\begin{equation}\label{sopr}
M^*(x,y) = \frac{xy}{M(x,y)},\qquad  x,y>0.
\end{equation}

Непосредственно проверяется, что сопряженное $M^*$ --- это также
среднее на множестве положительных чисел.

Известны достаточно широкие классы средних, удовлетворяющих
приведенным аксиомам. Перечислим основные из них.

1.  Степенные средние. Это самые известные средние:

\begin{equation}\label{M}
M(x,y) = M_{\alpha} (x,y) = \left( \frac{x^{\alpha}+y^{\alpha}}{2}
\right)^{\frac{1}{\alpha}} , -\infty \le \alpha \le \infty \ ,
\alpha \neq 0\ ;
\end{equation}
$$
M_{- \infty} (x,y) = \min (x, y)\ ,\  M_{0} = \sqrt{xy}\ , \
M_{\infty} (x, y) = \max (x, y)\ .
$$
Степенные средние образуют упорядоченную шкалу по параметру:
$$
\alpha_{1} > \alpha_{2} \Rightarrow M_{\alpha_1}(x,y) \ge
M_{\alpha_2} (x,y),\qquad \forall\,  x, y\ .
$$
Три исключительных значения $\{\alpha=-\infty, 0, +\infty\} $ могут
быть получены из неисключительных предельным переходом. Аксиомы
абстрактного среднего 1)--4) проверяются непосредственно.

Кроме симметричных средних  используются и несимметричные. Самые
известные из них --- это  средние арифметическое и геометрическое с
неотрицательными весами $\alpha,\beta$ :
$$
A_{\alpha,\beta}(x,y)=\alpha x + \beta y,\qquad
G_{\alpha,\beta}(x,y)=x^{\alpha} y^{\beta},\qquad \alpha +\beta=1.
$$
Неравенство между ними --- это  неравенство Янга, которое мы подробно рассмотрим в разделе \ref{kb}.  Аналогично
вводится весовое среднее степенное произвольного порядка.

Для сопряжённых этих средних получаем формулы
$$
A^*_{\alpha,\beta}(x,y)=\frac{1}{A_{\alpha,\beta}(\frac{1}{y},\frac{1}{x})},\qquad
G^*_{\alpha,\beta}(x,y)=G_{\beta,\alpha}(x,y).
$$

2. Средние Т.~Радо. Это средние следующего вида:
$$
R_{\beta}(x,y) = \left( \frac{x^{\beta + 1} - y^{\beta + 1}}{(\beta
+ 1)(x-y)}\right) ^{\frac{1}{\beta}},\qquad -\infty \le \beta \le
\infty ,\qquad  \beta \neq 0, -1;
$$
$$ R_{- \infty} (x,y) = \min (x, y),\qquad  R_{\infty} (x, y) = \max (x, y).
$$
Очевидно, что
$$
R_{-2} (x,y) = M_0 (x,y),\qquad R_{1} (x,y)=M_1 (x,y).
$$
Исключительные значения порождают пару гораздо менее известных средних:
логарифмическое
$$
R_{-1}(x,y)=L(x,y)=\frac{y-x}{\ln y-\ln x}
$$
и  "многоэтажное" \ среднее
$$
R_{0}(x,y)= \frac{1}{e} \left(
\frac{y^{y}}{x^{x}}\right)^{\frac{1}{y-x}}. \eqno(23)
$$

Средние Радо также образуют шкалу по параметру:
$$
\beta_{1} > \beta_{2} \Rightarrow R_{\beta_{1}} (x,y) \ge
R_{\beta_{2}} (x,y),\qquad  \forall x, y.
$$
Четыре исключительных значения  $\beta = \{-\infty, -1, 0,
+\infty\}$ могут быть получены из неисключительных предельным
переходом. Аксиомы абстрактного среднего  проверяются
непосредственно.

Средние Радо имеют простой аналитический смысл --- это промежуточные
значения в теореме Лагранжа о средних для логарифмической или
степенной функций. Например, неравенства для среднего
логарифмического из задачи \ref{zadL} (см. далее)  имеют смысл уточнений теоремы Лагранжа, из
которой для логарифмической функции следует только тривиальная
оценка:
$x\le L \le y$  \mbox{при} $x \le y$. Уже из оценки сверху через среднее арифметическое, полученной  В. Я.\,Буняковским как приложение его знаменитого неравенства,  следует неожиданный факт, что промежуточное значение из формулы Лагранжа для логарифмической функции может лежать не на всём указанном отрезке, как следует из стандартной формулировки, а только на левой половине отрезка, при этом оценка  задачи \ref{zadL} даёт дальнейшее существенное уточнение.

Приведённые неравенства для отдельных средних являются частными случаями общей теоремы, доказанной Тибором Радо, в которой устанавливаются неулучшаемые двусторонние оценки средних Радо через степенные средние и наоборот.

\begin{teo} \textbf{Теорема Радо}. Справедливы следующие двусторонние неулучшаемые
оценки средних Радо через степенные средние:
$$
M_{\frac{\alpha + 2}{3}}\le R_{\alpha}\le M_{0},\qquad  \mbox{при}\qquad
\alpha \in (-\infty, -2],
$$
$$
M_{0}\le R_{\alpha}\le M_{\frac{\alpha +2}{3}},\qquad  \mbox{при}\qquad
\alpha \in ~[-2, -1],
$$
$$
M_{\frac{\alpha\ln2}{\ln(1+\alpha)}}\le R_{\alpha} \le
M_{\frac{\alpha +2}{3}},\qquad \mbox{при}\qquad  \alpha \in (-1, -1/2],
$$
$$
M_{\frac{\alpha+2}{3}}\le R_{\alpha} \le M_{\frac{\alpha \ln
2}{\ln (1+\alpha)}},\qquad \mbox{при}\qquad \alpha \in ~[-1/2, 1),\\
$$
\mbox{(при $\alpha = 0$  последнее неравенство понимается в
предельном смысле}\\ $M_{\frac{2}{3}}\le R_0 \le M_{\ln 2}$),
$$M_{\frac{\alpha\ln
2}{\ln (1+\alpha)}}\le R_\alpha \le M_{\frac{\alpha +2}{3}},\qquad
\mbox{при}\qquad \alpha \in ~[1,\infty].\\
$$
\end{teo}

Набор средних степенных и Радо достаточно широк и обслуживает большинство прикладных задач. Тем не менее существуют средние, которые к ним не относятся. Примером может служить среднее Герона
$$
He(x,y)=\frac{x+\sqrt{xy}+y}{3}=\frac{1}{3}G(x,y)+\frac{2}{3}A(x,y)=
G_{\frac{1}{3},\frac{2}{3}}(G(x,y),A(x,y)).
$$
Среднее Герона является несколько неожиданным примером того, что несимметричное среднее от симметричных может опять оказаться симметричным средним.

3. Средние Джини и Лемера.  Итальянский статистик Коррадо Джини ввёл названные его именем средние с двумя параметрами  по формулам
$$
Gi_{u,v}(x,y)=\left( \frac{x^u+y^u}{x^v+y^v}\right)^{\frac{1}{u-v}},\qquad u\neq v,
$$
$$
Gi_{u,v}(x,y)=\exp \left(\frac{x^u\ln x+y^u \ln y}{x^u+y^u}\right),\qquad u=v\neq0,
$$
$$
Gi_{u,v}(x,y)=G(x,y),\qquad u=v=0.
$$

Важный частный случай средних Джини получается  при $v=u-1$.
Эти средние были переоткрыты  Д.~Лемером и имеют вид
$$
Le_u(x,y)=\frac{x^{u+1}+y^{u+1}}{x^u+y^u}.
$$

4. Квазиарифметические средние. Пусть дана неотрицательная монотонная функция $f(x)$, набор неотрицательных чисел $x=(x_1,x_2,\cdots,x_n)$ и неотрицательных весов
$p=(p_1,p_2,\cdots,p_n)$. \textit{Квазиарифметическим средним} называется выражение
$$
K_p(x)=f^{-1}\left( \sum_{k=1}^n p_k f(x_k)\right),\qquad \sum_{k=1}^n p_k=1.
$$
В частном случае $f(x)=x$ получаем обычное весовое среднее арифметическое, чем объясняется название.  Основные результаты для квазиарифметического среднего были получены в работах Колмогорова, Нагумо и Де Финетти. Они в существенном сводятся к доказательству двух основных свойств:\\

А1. Квазиарифметические средние совпадают тогда и только тогда, когда порождающие их функции связаны линейным соотношением.\\

А2. Непрерывное  квазиарифметическое среднее однородно тогда и только тогда, когда оно совпадает со средним степенным.\\

5. Итерационные средние.

Рассмотрим итерационный процесс при
заданных стартовых значениях $x_0, y_0$  и паре
абстрактных средних $(M, N)$:
$$
x_{n+1} = M(x_n, y_n), \ y_{n+1} = N(x_n, y_n).
$$

В общем случае получаем некоторую динамическую систему на плоскости,
ко\-то\-рая имеет интересное асим\-пто\-ти\-чес\-кое по\-ве\-де\-ние
(динамику) при $n\rightarrow+\infty$, но очень сло\-жна для изучения
даже при простейшем выборе пары средних $(M, N)$.

\textbf{Определение.} Пусть существует общий предел
последовательностей $x_n$ и $y_n$. Тогда он называется
\textit{итерационным средним} и обозначается
$$
\mu (M,N|\ x_0, y_0) = \mu(x_0, y_0) = \lim x_n = \lim y_n.
$$

Непосредственно проверяется, что итерационное среднее
пары абстрактных средних также является новым абстрактным средним и
наследует соответствующие аксиомы 1)--4). Также принято другое обозначение вместо $\mu(M,N)$, а именно $M\otimes N$.

Самое известным итерационным средним является \\арифметико--геометрическое (AGM = arithmetic--geometric mean), изученное Лежандром и Гауссом. Оно получается при выборе $M=M_1$, $N=M_0$ \   и выражается по формуле
$$
\mu(M,N|\ x_0,y_0) = \frac{\frac{\pi}{2}\ x_0}{K\left(\sqrt{1-\left(
\frac{y_0}{x_0}\right)^{2}}\right)},\qquad 0 < y_0 <x_0,
$$
где $K(x)$ --- полный эллиптический интеграл Лежандра первого рода, см. \cite{Fikh} .

Таким образом, можно сделать основной  вывод, что
\textit{существует значительное число конкретных примеров
абстрактных средних}.

\subsection{Неравенства для средних в комплексной плоскости}

Одна из основных структурных теорем алгебры запрещает сравнивать комплексные числа, так как утверждает существование единственного с точностью до изоморфизма упорядоченного числового поля, а это место уже занято действительными числами. Остаётся переходить к модулям.

Итак, пусть выполнено неравенство $f(x,y) \le g(x,y)$ для неотрицательных функций $f(x,y) \ge 0, g(x,y)\ge 0$ от \textit{действительных} переменных $x,y \ge 0$. Назовём  \textit{комплексификацией действительного неравенства} набор следующих неравенств \eqref{com1}--\eqref{com2} для \textit{комплексных} переменных $z,w\in \mathbb{C}$:
\begin{equation}\label{com1}
|f(z,w)| \le |g(z,w)|,
\end{equation}
\begin{equation}\label{com2}
|f(z,w)| \le g(|z|,|w|),
\end{equation}
\begin{equation}\label{com3}
f(|z|,|w|) \le |g(z,w)|,
\end{equation}
\begin{equation*}
f(|z|,|w|) \le g(|z|,|w|).
\end{equation*}

Рассматривать последнее из приведённых выше неравенств не имеет смысла, так как оно совпадает с исходным. На самом деле выписан далеко не полный комплект, что отражают многоточия, ведь модули можно вставить в трёх местах справа и трёх местах слева, итого получается $8 \times 8=64$ варианта, правда, некоторые из них отпадут. Например, неравенство $f(|z|,|w|) \le |g(|z|,w)|$ может быть добавлено в список, а вариант $f(z,|w|) \le |g(|z|,w)|$ следует забраковать, так как слева стоит комплексное число. Отдельно нужно рассматривать, когда неравенства переходят в равенства.
Мы ограничимся комплектом для комплексификации \eqref{com1}--\eqref{com3}.

Разумеется, приведённая схема подходит и для неравенств от одной переменной, и от нескольких, а не только от двух.

Теперь можно применять введённую комплексификацию ко многим известным неравенствам. Мы рассмотрим самое простое --- неравенство между средним арифметическим и геометрическим. Результаты получаются достаточно неожиданные и красивые. В принципе, процесс можно продолжить на две оставшиеся по теореме Гурвица нормированные гиперкомплексные алгебры --- кватернионы и октавы.

Итак, дано исходное неравенство
$$
f(x,y)=\sqrt{xy} \le g(x,y)=(\frac{x+y}{2}), x\ge 0, y\ge 0.
$$
Корень нам не совсем подходит, поэтому, чтобы избежать многозначных функций, возведём это неравенство в квадрат. Его комплексификации выписываются тогда при $z,w\in \mathbb{C}$ так:
\begin{equation}\label{com4}
|zw| \le \left|\frac{z+w}{2}\right|^2,\qquad
|zw| \le {\left(\frac{|z|+|w|}{2}\right)}^2,\qquad
|zw| \le \left|\frac{|z|+w}{2}\right|^2.
\end{equation}
Для рассматриваемого неравенства весь полный боекомплект из 64 вариантов сводится к трём \eqref{com4}.

Мы хотим рассчитать и потом изобразить на графиках те множества, для которых выполняются неравенства \eqref{com4}. Сделать это в 4--мерном пространстве сложновато, поэтому упростим задачу и используем однородность в \eqref{com4}. Тогда для новой комплексной переменной $s=\frac{w}{z}$ получим
\begin{equation}\label{com5}
|s| \le \left|\frac{s+1}{2}\right|^2,\qquad
|s| \le {\left(\frac{|s|+1}{2}\right)}^2,\qquad
|s| \le \left|\frac{\frac{w}{|z|}+1}{2}\right|^2.
\end{equation}
Два первых неравенства теперь зависят от одной переменной, их мы и продолжим рассматривать. Что делать с третьим, которое опять зависит от двух комплексных переменных, непонятно, и мы его рассматривать не будем.

Исследуем первое из неравенств \eqref{com5}. После подстановки $s=x+iy$ и некоторых упрощений, получаем, что равенство в нём достигается на некоторой кривой 4 порядка с уравнением
\begin{equation}\label{com6}
x^4+y^4+2x^2y^2+4x^3+4xy^2-10x^2-14y^2+4x+1=0.
\end{equation}

C помощью Maple и некоторых расчётов можно сформулировать выводы.

\begin{teo}
\begin{enumerate}
\item  На комплексной плоскости существует замкнутая ограниченная кривая, на которой в первом из неравенств \eqref{com5} достигается равенство, вне её выполнено неравенство, а в ограниченной области внутри кривой справедливо обратное неравенство. Указанная кривая определяется уравнением \eqref{com6}.\\

\item  Первое из неравенств \eqref{com5} для действительных $s$ справедливо не только для $s{>}0$, как следовало бы ожидать, а при $s\in[-\infty,-3-2\sqrt2]\cup[-3+2\sqrt2,\infty]$.\\

\item  Первое из неравенств \eqref{com5} для чисто мнимых $s$ выполняется при ${\rm Im}(s)\in [-\infty,-2-\sqrt3]\cup[-2+\sqrt3,2-\sqrt3]\cup[2+\sqrt3,\infty]$.\\

\item  Уравнения внешней и внутренней кривых в полярных координатах имеют вид $r=2-\cos{\phi} \pm \sqrt{(2-\cos{\phi})^2-1}$.
\end{enumerate}
\end{teo}

\begin{figure}[h!]
	\begin{minipage}[h]{0.49\linewidth}
		\center{\includegraphics[scale=0.25]{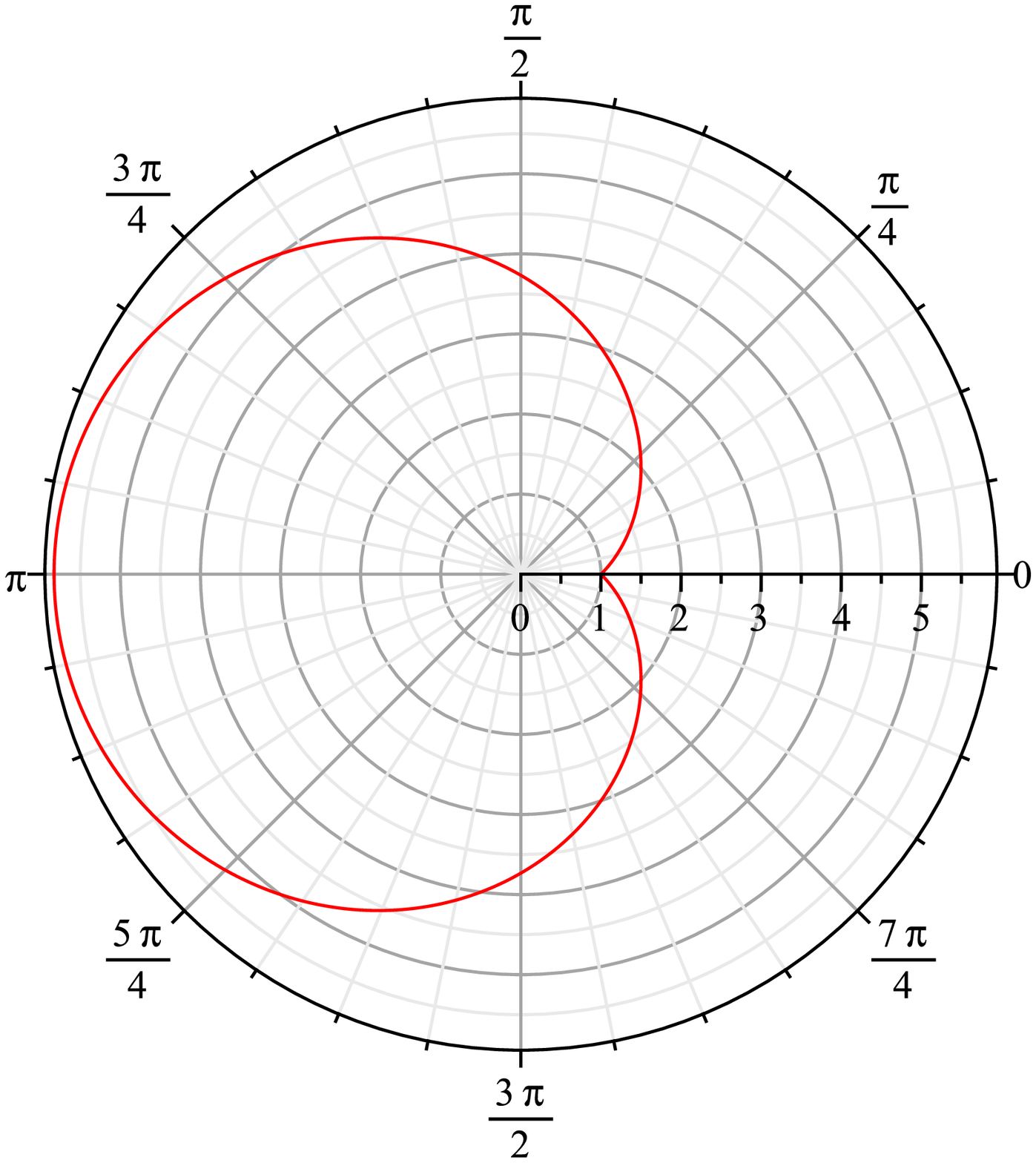} \\ a)}
	\end{minipage}
	\hfill
	\begin{minipage}[h!]{0.49\linewidth}
		\center{\includegraphics[scale=0.25]{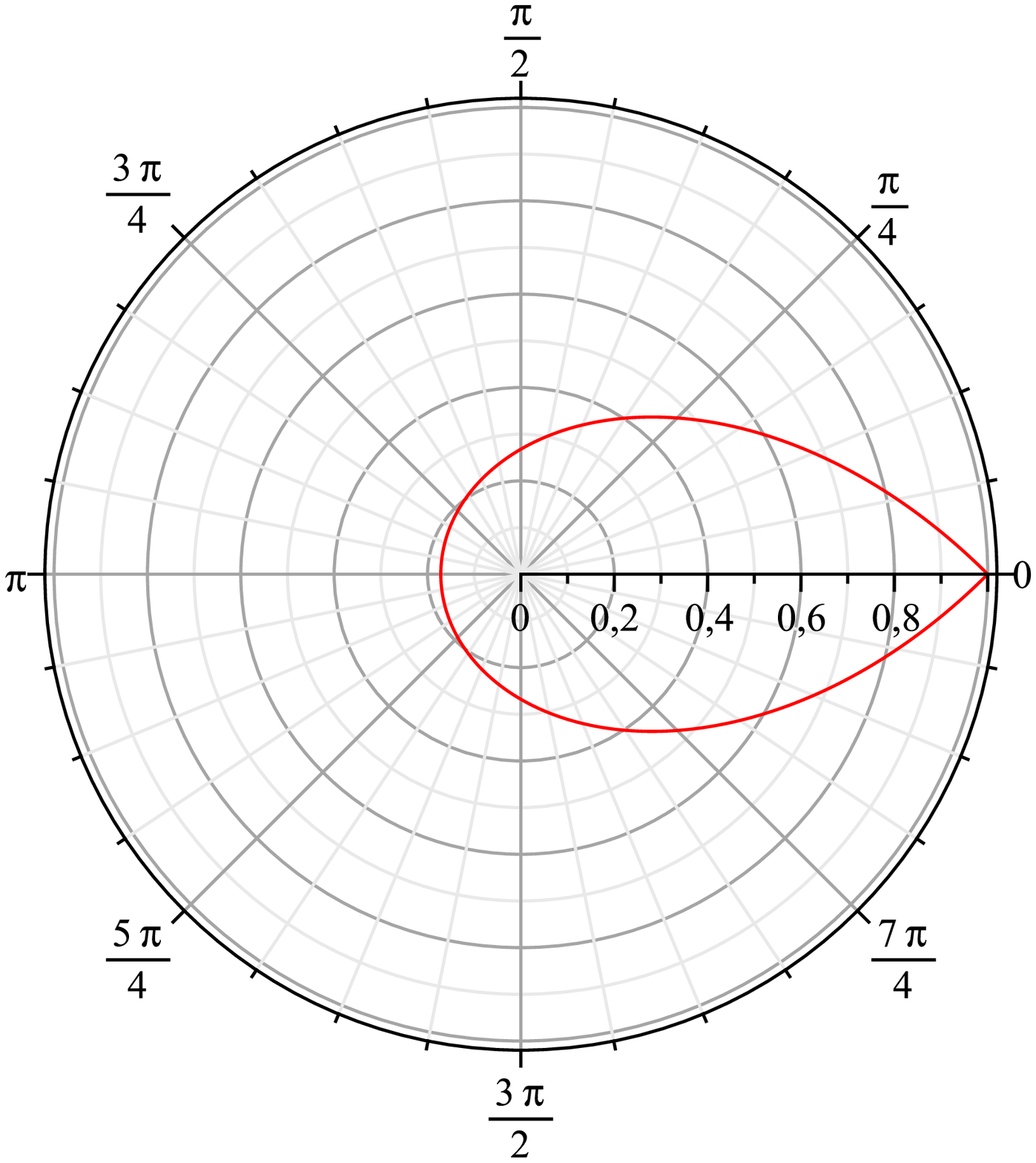} \\ b)}
	\end{minipage}
	\caption{Первая a) и вторая b) части кривой $r=2-\cos{\phi} \pm \sqrt{(2-\cos{\phi})^2-1}$ в полярных координатах.}
	\label{ris1}
\end{figure}

\begin{figure}[h!]
\begin{center}
\includegraphics[scale=0.3]{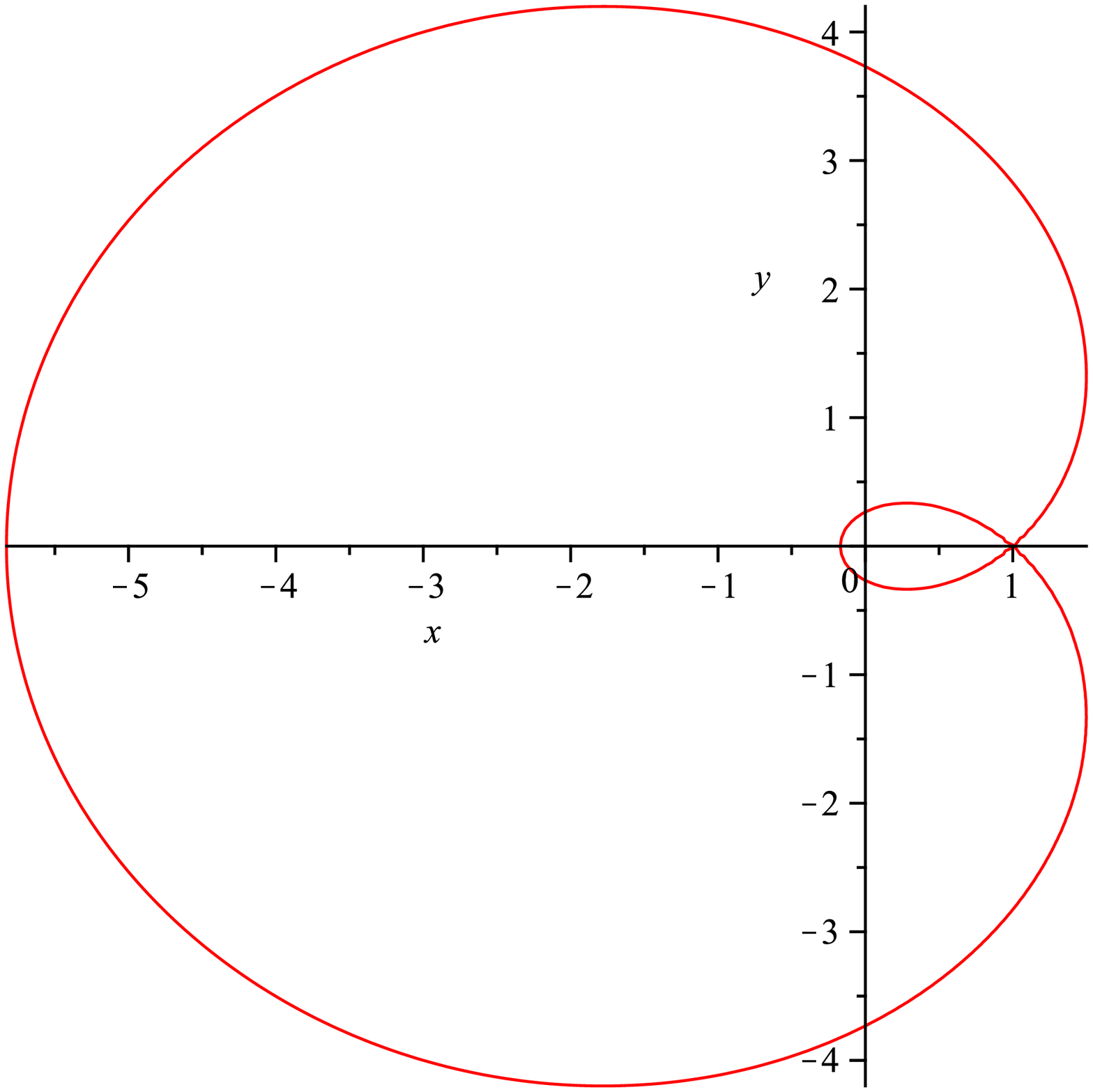}
\caption{Полный график  разделяющей кривой   в прямоугольных координатах.}\label{ris2}
\end{center}
\end{figure}

Приведём численные значения $-3-2\sqrt2\approx-5,8284, -3+2\sqrt2\approx-0,1716, -2-\sqrt3\approx-3,7321, -2+\sqrt3\approx-0,2679,  2-\sqrt3\approx0,2679, 2+\sqrt3\approx3,7321.$
Графики $r=2-\cos{\phi} \pm \sqrt{(2-\cos{\phi})^2-1}$ сначала отдельно двух кусков в полярных координатах, а затем всей кривой в прямоугольных координатах приведены  на Рис. \ref{ris1} и \ref{ris2}.
Получившиеся кривые похожи на некоторые известные кривые 4 порядка,
а именно, улитки Паскаля, кардиоиды и овалы Декарта, но к ним не сводятся.

Что же касается второго из неравенств \eqref{com5}, то с ним всё ясно. Это неравенство только для одного радиуса и оно выполняется на границе и вне единичного круга и не выполняется внутри его.

Таким образом, указанная процедура комплексификации неравенств может рассматриваться как источник новых плоских кривых. Их можно назвать их разделяющими для комплексификации соответствующего неравенства.

Следует сказать, что неравенства для средних в комплексной плоскости --- это не только интересная забава.  В знаменитой работе В.~Бекнера, в которой получены точные значения норм преобразования Фурье в
шкале пространств $L_p$, используются две основные конструкции:
дробное (или квадратичное) преобразование Фурье и некоторое неравенство для средних степенных  в комплексной
плоскости. Точную константу, которую получил Бекнер, искали и многие до него. Наибольших успехов добился наш математик К.~И.~Бабенко, решивший задачу в частном случае $p=2^n$. Ему принадлежит идея применять в доказательстве специальное квадратичное преобразование Фурье. Вторая идея принадлежит Бекнеру, он свёл оставшуюся часть задачи к некоторому неравенству между степенными средними в комплексной плоскости. Нужно было доказать, что это неравенство выполняется на определённом отрезке мнимой оси, что Бекнеру и удалось сделать.

\subsection{Неравенства для средних значений}

\begin{zd}
Доказать неравенства между основными степенными средними:
$$
M_{-1}(x,y)=H(x,y) \le M_0(x,y)=G(x,y)  \le M_1(x,y)=
$$
$$
=A(x,y) \le M_2(x,y)=Q(x,y).
$$
\end{zd}

\begin{zd}Доказать формулу для сопряжённого степенного среднего
$$
(M_\alpha)^*=M_{-\alpha}.
$$
\end{zd}

\begin{zd}
Доказать монотонность степенных средних по параметру.
\end{zd}

\begin{zd}\hspace{-2mm}*\label{zadL}
Доказать неравенства для некоторых  средних Радо, логарифмического и "многоэтажного"\ :
$$
M_{0}(x,y) = \sqrt{xy}\le L(x,y) = \frac{x-y}{\ln x - \ln y}
\le M_{\frac{1}{3}}(x,y)=\left(
\frac{x^{\frac{1}{3}}+y^{\frac{1}{3}}}{2}\right)^{3},
$$
$$
M_{\frac{2}{3}}(x,y) = \left(
\frac{x^{\frac{2}{3}}+y^{\frac{2}{3}}}{2}\right)^{\frac{2}{3}} \le
R_{0}(x,y) =
$$
$$
 = \frac{1}{e}\left(
\frac{y^{y}}{x^{x}}\right)^{\frac{1}{y-x}} \le M_{\ln2}(x,y) =
\left( \frac{{x^{\ln2}}+y^{\ln2}}{2}\right)^{\frac{1}{\ln2}}.
$$
\end{zd}

\begin{zd}\hspace{-2mm}***
Доказать, что в неравенствах предыдущей задачи  порядки средних степенных $0, 1/3$ и $2/3, \ln2$ неулучшаемы.
\end{zd}

\begin{zd}Доказать формулу для сопряжённого  среднего Радо
$$
R^*_{\beta}(x,y)=\frac{1}{R_{\beta}\left( \frac{1}{x}, \frac{1}{y}\right)}.
$$
\end{zd}

\begin{zd}
Доказать монотонность  средних Радо по параметру.
\end{zd}

\begin{zd}\hspace{-2mm}***
Докажите, что только пять классических средних входят одновременно в оба множества степенных средних и Радо:
$$
M_{-\infty} = R_{-\infty}, \quad M_{0} = R_{-2},\quad M_{\frac{1}{2}} =
R_{\frac{1}{2}},\quad M_{1} = R_{1},\quad M_{\infty} = R_{\infty}.
$$
\end{zd}

\begin{zd}\hspace{-2mm}**
Доказать неравенства для арифметико--геометрического среднего Гаусса
$$
L(x, y) = R_{-1}(x,y)\le\mu(x,y)\le
R_{\frac{1}{2}}(x,y)=M_{\frac{1}{2}}(x,y).
$$
\end{zd}

Итак, средних очень много. Тем удивительнее, что можно получить описание всех "хороших"\ средних от двух переменных.

\begin{zd}\hspace{-2mm}***
Доказать, что произвольное абстрактное среднее, обладающее свойствами несмещённости, однородности, симметричности,
монотонности и непрерывности по обоим аргументам представляется в виде
\begin{equation}\label{sr}
M (x,y)= (x+y)h\left(\ln \frac{y}{x}\right),
\end{equation}
 где  $h(t)$ является определенной на всей оси непрерывной
 четной функцией, удовлетворяющей следующим условиям при $t\ge
0$:
$$
h(0)=\frac{1}{2} ,\\
$$
$$
\frac{e^{t_1}(e^{t_2}+1)}{e^{t_2}(e^{t_1}+1)}\ \le \
\frac{h(t_1)}{h(t_2)} \le  \frac{(e^{t_2}+1)}{(e^{t_1}+1)}, \qquad
t_1\le t_2.
$$
Справедливо и обратное: каждой функции $h(t)$ с указанным
набором свойств соответствует по формуле \eqref{sr} несмещённое,
однородное, симметричное, монотонное и непрерывное по обоим
аргументам среднее.
\end{zd}

\begin{zd}\hspace{-2mm}***
 Исследовать предложенным выше методом комплексификации неравенство $\sin x \le x,$ $x\ge 0$.
\end{zd}

\section{Неравенства Коши--Буняковского, Минковского, Янга, Гёльдера и некоторые их уточнения}\label{kb}

Приведём классические неравенства Коши--Буняковского и Минковского для конечных сумм, см., например, \cite{Kud}.

\begin{teo}
Неравенство Коши--Буняковского имеет вид
\begin{equation}\label{CB}
  \left|\sum _{i=1}^{n}u_{i}v_{i}\right| \leq \sqrt{\left(\sum _{i=1}^{n}u_{i}^{2}\right)\left(\sum _{i=1}^{n}v_{i}^{2}\right)}.
\end{equation}
\end{teo}
{\it Доказательство:} Действительно, если все $u_i=0$, $i=1,2,...,n$, то неравенство \eqref{CB} очевидно, так как обе его части обращаются в нуль. Если же
$\sum\limits_{i=1}^{n}u_{i}^{2}>0$, то рассмотрим функцию
\begin{equation}\label{KT}
F(t)=\sum\limits_{i=1}^{n}(u_{i}t+v_{i})^2=t^2 \sum\limits_{i=1}^{n}u_{i}^{2}+2t\sum\limits_{i=1}^{n}u_{i}v_i+\sum\limits_{i=1}^{n}v_{i}^2.
\end{equation}
Очевидно, что $F(t)\geq0$, тогда дискриминант квадратного трехчлена  \eqref{KT} либо равен нулю, либо отрицателен:
$$
\left(\sum _{i=1}^{n}u_{i}v_{i}\right)^{2}-\left(\sum _{i=1}^{n}u_{i}^{2}\right)\left(\sum _{i=1}^{n}v_{i}^{2}\right)\leq 0.
$$
Перенеся второе слагаемое в правую часть и извлекая квадратный корень, получим \eqref{CB}.

\vskip 0.5 cm

\begin{teo}
	Неравенство Минковского имеет вид
	\begin{equation}\label{Min}
	\sqrt{\sum\limits_{i=1}^{n}(u_{i}+v_{i})^2}\leq \sqrt{\sum _{i=1}^{n}u_{i}^{2}}+\sqrt{\sum _{i=1}^{n}v_{i}^{2}}.
	\end{equation}
\end{teo}
{\it Доказательство:} Оценим сумму
$$
\sum\limits_{i=1}^{n}(u_{i}+v_{i})^2=\sum\limits_{i=1}^{n}u_i^2+2\sum\limits_{i=1}^{n}u_iv_i+\sum\limits_{i=1}^{n}u_i^2.
$$
Применяя неравенство \eqref{CB}, получим
$$
\sum\limits_{i=1}^{n}(u_{i}+v_{i})^2\leq\sum\limits_{i=1}^{n}u_i^2+2\sqrt{\sum _{i=1}^{n}u_{i}^{2}}\sqrt{\sum _{i=1}^{n}v_{i}^{2}}+\sum\limits_{i=1}^{n}u_i^2=
$$
$$
=\left(\sqrt{\sum _{i=1}^{n}u_{i}^{2}}+\sqrt{\sum _{i=1}^{n}v_{i}^{2}}\right)^2.
$$
Извлекая из обеих частей квадратный корень, получим \eqref{Min}.

Следующее неравенство доказано Вильямом Янгом, в русской литературе оно часто называется неравенством Юнга.

\begin{teo} \textbf{Неравенство Янга.} Справедливо неравенство при указанных условиях
\begin{equation}\label{Yng1}
xy\le \frac{x^{p} }{p} +\frac{y^{q} }{q},\\
x\ge 0, y\ge 0,     p>1,   \frac{1}{p} +\frac{1}{q} =1.
\end{equation}
\end{teo}

\begin{zd} Доказать неравенство Янга.
\end{zd}

Неравенство Янга может быть сформулировано с использованием средних значений. Тогда оно превращается в другое классическое неравенство: весовое среднее геометрическое не превосходит весового среднего арифметического:
 $$
 u^{a} v^{b} \le au+bv,\begin{array}{cc} {} & {} \end{array}a\ge 0,b\ge 0,a+b=1,u\ge 0,v\ge 0.
 $$

Далее в этом пункте излагается материал, следуя обзору \cite{Sit2}, а также работам \cite{Sit3--Sit7}.

Принято считать, что круг вопросов, связанных с простейшим неравенством \eqref{Yng1}, исследован с исчерпывающей полнотой. Но это не совсем так. Если обратить внимание на несимметричную форму правой части \eqref{Yng1}, то становится ясно, что \textit{неравенство Янга --- это не одно, а пара неравенств}, причём второе из них, пропущенное  в \eqref{Yng1}, должно иметь такой вид:
\begin{equation}\label{Yng2}
xy\le \frac{x^{q} }{q} +\frac{y^{p} }{p}.
\end{equation}

Тогда  возникает задача о сравнении неравенств \eqref{Yng1} и \eqref{Yng2}, то есть об отыскании минимума их правых частей. Приведем численные примеры, которые показывают, что действительно между правыми частями  \eqref{Yng1} и \eqref{Yng2} может существовать достаточно большой разброс.

\begin{pr}  Выберем значения $x=5, y=130, p=4, q=4/3$; тогда
$$
xy = 650,\qquad    \frac{x^{p} }{p} +\frac{y^{q} }{q} \approx 650,16502,$$
$$\frac{x^{q} }{q} +\frac{y^{p} }{p} \approx 71402508.$$

\end{pr}
В этом случае неравенство \eqref{Yng1} лучше (на пять порядков!).

\begin{pr}  Выберем значения $x=0,2, y=0,5, p=4, q=4/3$; тогда\\
$$xy = 0,1,\qquad   \frac{x^{p} }{p} +\frac{y^{q} }{q} \approx 0,29803,$$
$$\frac{x^{q} }{q} +\frac{y^{p} }{p} \approx 0,10334.$$

\end{pr}
А в этом случае неравенство \eqref{Yng2} лучше (примерно в три раза).

\begin{zd}
Проверьте приведённые результаты вычислений.
\end{zd}

Два приведённых примера иллюстрируют два типичных случая (хотя есть и третий!). Мы сразу сформулируем общий результат, см. \cite{Sit2}. Трудность была только в том, чтобы догадаться и найти данную оценку, после этого доказательство превращается в не очень сложное упражнение на применение производной. Мы предлагаем читателю выполнить это небольшое исследование самостоятельно, для его удобства мы разобьём решение ниже на несколько задач .

Без ограничения общности далее будем предполагать, что выполнены условия  $x\ge 0,$ $ y\ge 0$,     $p>1$,   $\frac{1}{p} +\frac{1}{q} =1.$

\begin{zd} Доказать, что если  $y\ge x\ge 1,$  то оценка \eqref{Yng1} лучше, чем \eqref{Yng2}, то есть выполнены неравенства
\begin{equation}\label{Yng3}
xy\le \frac{[\min (x,y)]^{p} }{p} +\frac{[\max (x,y)]^{q} }{q} \le \frac{[\max (x,y)]^{p} }{p} +\frac{[\min (x,y)]^{q} }{q}.
\end{equation}
\end{zd}

\begin{zd}Доказать, что если   $1\ge y\ge x\ge 0,$  то оценка \eqref{Yng2} лучше, чем \eqref{Yng1}, то есть выполнены неравенства
\begin{equation}\label{Yng4}
xy\le \frac{[\max (x,y)]^{p} }{p} +\frac{[\min (x,y)]^{q} }{q} \le \frac{[\min (x,y)]^{p} }{p} +\frac{[\max (x,y)]^{q} }{q}.
\end{equation}
\end{zd}

\begin{zd}\hspace{-2mm}***  Доказать, что если   $y\ge 1\ge x\ge 0,$  то при данном $x$ существует единственное критическое значение $y=y_{cr}$, которое является решением трансцендентного уравнения
\begin{equation}\label{Yng5}
\frac{x^{p} }{p} -\frac{x^{q} }{q} =\frac{y^{p} }{p} -\frac{y^{q} }{q}.
\end{equation}
В этом случае при  $1\le y\le y_{cr}$  оценка \eqref{Yng2} лучше, чем \eqref{Yng1}, то есть выполнены неравенства \eqref{Yng4}, а при  $y\ge y_{cr}$  оценка \eqref{Yng1} лучше, чем \eqref{Yng2}, то есть выполнены неравенства \eqref{Yng3}.
\end{zd}

Иначе говоря, если два числа лежат с одной стороны от единицы, то лучше одно из неравенств, с другой стороны---лучше другое, а если единица разделяет два числа, то реализуются оба случая.

Теперь приведем обещанный численный пример на пропущенный третий случай.

\begin{pr} Выберем значения  $x=0,5,$ $p=4$, $q=4/3$. Тогда расчет  дает критическое значение $y_{cr} \approx 1,35485.$
Выберем  $x=0,5$, $y=1$, $3<y_{cr}$; тогда получаем $xy =0,65$,
\begin{equation}
\frac{x^{p} }{p} +\frac{y^{q} }{q} \approx 1,07973, \frac{x^{q} }{q} +\frac{y^{p} }{p} \approx 1,01166.
\end{equation}
Как и следует из вышеизложенного, в этом случае лучше оценка \eqref{Yng2}.
\end{pr}
\begin{pr} Выберем значения  $x=0,5,\  y=1,4> y_{cr}$;  тогда получаем $xy=0,7$,
\begin{equation}
\frac{x^{p} }{p} +\frac{y^{q} }{q} \approx 1.19025,\frac{x^{q} }{q} +\frac{y^{p} }{p} \approx 1.25804.
\end{equation}
Как и следует из вышеизложенного, в этом случае лучше оценка \eqref{Yng1}.
\end{pr}

\begin{zd}\hspace{-2mm}***
Разработать метод численного нахождения критического значения $y_{cr}$ с использованием метода Ньютона.
\end{zd}

\begin{zd}
Вывести из неравенства Янга при условии $p{>}1$, $\frac{1}{p}{+}\frac{1}{q}{=}1$ неравенство Гёльдера
\begin{equation}\label{Hol}
  \left|\sum _{i=1}^{n}u_{i}v_{i}\right|
   \leq
   \left(\sum _{i=1}^{n}u_{i}^{p}\right)^{1/p}
  \left(\sum _{i=1}^{n}v_{i}^{q}\right)^{1/q}.
\end{equation}
\end{zd}

На самом деле то, что мы называем неравенством Гёльдера, было первоначально доказано замечательным английским математиком Л.\,Роджерсом, чуть позже О.\,Гёльдером; при этом Гёльдер цитировал результат Роджерса, и на самом деле они оба доказали не тот вариант неравенства \eqref{Hol}, к которому мы привыкли. А общепризнанный сейчас вариант неравенства \eqref{Hol} был впервые сформулирован и доказан Ф.\,Риссом. Поэтому правильный вариант названия --- неравенство Роджерса--Гёльдера--Рисса.

\newpage

\section{ Ссылки}

\end{document}